\documentclass[a4paper,fleqn]{scrartcl}
\usepackage{times,amsmath}
\usepackage{cite}
\usepackage{ calc, color}
\usepackage[]{graphicx}
\usepackage[affil-it]{authblk}
\begin{document}
\title{Modelling and simulation of acrylic bone cement injection and curing within the framework of vertebroplasty
       \footnote{This is the pre-peer reviewed version of the following article: 
    R. Landgraf, J. Ihlemann, S. Kolmeder, A. Lion, H. Lebsack, C. Kober. Modelling and simulation of acrylic bone cement injection and curing within the framework of vertebroplasty. ZAMM - Zeitschrift f\"ur angewandte Mathematik und Mechanik, DOI 10.1002/zamm.201400064 (2015), 
	which has been published in final form at \textcolor{blue}{http://onlinelibrary.wiley.com/doi/10.1002/zamm.201400064/abstract}. }}

\author{Ralf Landgraf$^1$%
            \thanks{Corresponding author\quad 
						E-mail:~\textsf{ralf.landgraf@mb.tu-chemnitz.de}}\,\,}
\author{J\"orn Ihlemann$^1$}
\author{Sebastian Kolmeder$^2$}
\author{Alexander Lion$^2$}
\author{Helena Lebsack$^3$}
\author{Cornelia Kober$^3$}

\affil{	\small $^1$Department of Solid Mechanics, Faculty of
		Mechanical Engineering, Chemnitz University of  
		Technology, Chemnitz, Germany; 
		$^2$ Institute of Mechanics, Faculty of Aerospace 
		Engineering, University of the Federal Armed Forces of 
		Germany, Neubiberg, Germany;
		$^3$ Faculty of Life Sciences, Hamburg University of Applied 
		Sciences, Hamburg, Germany}
\date{March 27, 2014}
\maketitle                   
\vspace{-5ex}
\renewenvironment{abstract}{%
   \begin{minipage}{0.95\textwidth}}
   {\par\noindent
   \end{minipage}
}
\begin{abstract}
The minimal invasive procedure of vertebroplasty is a surgical technique to treat compression fractures of vertebral bodies. During the treatment, liquid bone cement gets injected into the affected vertebral body and therein cures to a solid. In order to investigate the treatment and the impact of injected bone cement, an integrated modelling and simulation framework has been developed. The framework includes (i) the generation of microstructural computer models based on microCT images of human cancellous bone, (ii) computational fluid dynamics (CFD) simulations of bone cement injection into the trabecular structure and (iii) non-linear finite element (FE) simulations of the subsequent bone cement curing. A detailed description of the material behaviour of acrylic bone cements is provided for both simulation stages. A non-linear process-depending fluid flow model is chosen to represent the bone cement behaviour during injection. The bone cements phase change from a highly viscous fluid to a solid is described by a non-linear viscoelastic material model with curing dependent properties. To take into account the distinctive temperature dependence of acrylic bone cements, both material models are formulated in a thermo-mechanically coupled manner. Moreover, the corresponding microstructural CFD- and FE-simulations are performed using thermo-mechanically coupled solvers. An application of the presented modelling and simulation framework to a sample of human cancellous bone demonstrates the capabilities of the presented approach.
\end{abstract}

\def\etal{\mbox{et al.}}
\newbox\JIGAMMa
\newbox\JIGAMMb
\newbox\JIGAMMc
\newbox\JIGAMMd
\newbox\JIGAMMz
\newbox\LLbox
\newbox\LLboxh
\newbox\SLhilfbox
\newbox\SLubox
\newbox\SLobox
\newbox\SLergebnis
\newbox\TENbox
\newif\ifSLoben
\newif\ifSLunten
\newdimen\JIGAMMdimen
\newdimen\JIhsize\relax\JIhsize=\hsize
\newdimen\SLrandausgleich
\newdimen\SLhoehe
\newdimen\SLeffbreite
\newdimen\SLuvorschub
\newdimen\SLmvorschub
\newdimen\SLovorschub
\newdimen\SLsp
\def\rhotilzurho{{\JI\frac{\STAPEL\varrho!^\SLtilde\!}{\JI\!\varrho}}}
\def\rhozurhotil{{\JI\varrho\over\JI\STAPEL\varrho!^\SLtilde}}
\def\pkt{\cdot}
\def\ppkt{\mathbin{\mathord{\cdot}\mathord{\cdot}}}
\setbox\JIGAMMa = \hbox{$\scriptscriptstyle c$} \setbox\JIGAMMz =
\hbox{\hskip-.35pt\vrule width .25pt\hskip-.35pt
                        \vbox to1.2\ht\JIGAMMa{\vskip-.125pt
                             \hrule width1.2\ht\JIGAMMa height.25pt
                             \vfill
                             \hrule width1.2\ht\JIGAMMa height.25pt
                             \vskip-.125pt}%
                        \hskip-.125pt\vrule width .25pt\hskip-.125pt}
\def\Oldroy#1#2#3{\STAPEL{#1}!_\SLstrich!_\SLstrich!^\circ{}
                  \ifx #2,{}_{\copy\JIGAMMz}%
                  \else \mskip1mu{}^{\copy\JIGAMMz}\fi
                  \mskip1mu\ifx #3,{}_{\copy\JIGAMMz}%
                           \else {}^{\copy\JIGAMMz}\fi }
\def\OP#1#2{\ifnum#1=1{\rm S}
            \else\ifnum#1=2{\rm S}^\JIv
                 \else\ifnum#1=3{\rm S}^T
                      \else{{\rm S}^T}^\JIv
            \fi\fi\fi\LL{{#2}}\RR}
\edef\JIminus{{\setbox\JIGAMMa=\hbox{$\scriptstyle x$}%
           \hbox{\hskip .10\wd\JIGAMMa
                 \vbox{\hrule width .6\wd\JIGAMMa height .07\wd\JIGAMMa
                       \vskip.53\ht\JIGAMMa}%
                 \hskip .10\wd\JIGAMMa}}}
\edef\JIv{{\JIminus 1}}
\def\JI{\displaystyle}
\def\JIha{{1\over 2}}
\def\JIfolgt{\quad\Rightarrow\quad}
\def\LL#1\RR{\setbox\LLbox =\hbox{\mathsurround=0pt$\displaystyle
                                              \left(#1\right)$}%
       \setbox\LLboxh=\hbox{\mathsurround=0pt%
                  $\displaystyle{\left(%
                      \vrule width 0pt height\ht\LLbox depth\dp\LLbox
                      \right)}$}%
       \left(\hskip-.3\wd\LLboxh\relax\copy\LLbox
              \hskip-.3\wd\LLboxh\relax\right)}
\def\ZBOX#1#2#3{\def#3{}%
                \setbox#1 = #2
                \def#3{ to \wd#1}%
                \setbox#1 = #2}
\def\SLdreieck{\setbox\TENbox=\hbox{\fontscsy\char 52}
                 \dp\TENbox = 0pt
                 \hbox{\hskip -2\SLrandausgleich
                       \box\TENbox
                       \hskip -2\SLrandausgleich}}
\def\SLtilde{\setbox\TENbox=\hbox{\fontscex\char 101}
                    \vbox{\vskip-.03\ht\TENbox
                          \hbox{\hskip -1\SLrandausgleich
                                \copy\TENbox
                                \hskip -1\SLrandausgleich}
                          \vskip -.86\ht\TENbox}}
\def\SLstrich{\vrule width \SLeffbreite height.4pt}
\def\SLpunkt{{\vbox{\hbox{$\displaystyle.$}\vskip.03cm}}}
\def\SLabstand{\vskip .404pt}
\def\SLzwischen{\vskip 1.372pt}
\font\fontscsy=cmsy6 \font\fontscex=cmex10 scaled 1200
\def\STAPEL#1{\def\SLkern{#1}%
              \futurelet\next\SLpruef
               A_0   _0    :B_0   _-.17 :C_.05 _-.15 :D_0   _-.2
              :E_0   _-.2  :F_0   _-.21 :G_0   _-.15 :H_0   _-.23
              :I_.2  _.15  :J_.05 _-.1  :K_0   _-.22 :L_0   _-.1
              :M_0   _-.23 :N_0   _-.25 :O_.05 _-.2  :P_0   _-.21
              :Q_.05 _-.2  :R_0   _-.03 :S_0   _-.15 :T_.2  _0
              :U_.1  _-.1  :V_.1  _-.15 :W_.1  _-.2  :X_0   _-.22
              :Y_.16 _-.15 :Z_0   _-.25
              :a_.05 _-.05 :b_.05 _0    :c_.05 _.05  :d_0   _-.05
              :e_.07 _0    :f_0   _-.15 :g_.04 _-.2  :h_0   _-.07
              :i_.05 _0    :j_.08 _-.1  :k_0   _-.1  :l_.2  _.15
              :m_0   _-.1  :n_0   _-.1  :o_0   _-.1  :p_.15 _0
              :q_.1  _0    :r_.1  _-.1  :s_0   _-.2  :t_.1  _.05
              :u_0   _-.1  :v_0   _-.2  :w_0   _-.2  :x_.04 _-.14
              :y_.15 _-.05 :z_0   _-.15
              :\mit\Phi_.08 _-.1    :\mit\Omega_0 _-.2   :\varXi_.00 _-.2
              :\alpha_0 _-.2        :\gamma_.1 _-.1      :\varepsilon_.1 _-.1
              :\epsilon_.05 _-.05   :\eta_.05 _-.15      :\lambda_0 _0
              :\mu_0 _-.25          :\nu_0 _-.2          :\varSigma_-.03 _-.2
              :\varrho_.00 _-.2     :\sigma_.1 _-.2      :\tau_.15 _-.15
              :\varphi_.2 _-.1      :\omega_.1 _-.1      :\mit\Gamma_-.1 _-.1
              :\Lambda_0 _0         :\Gam_0 _0           :\Lam_0 _0
              :\SLsuchende
              \def\SLtrick{\noexpand\SLtrick\noexpand}%
                \def\SLdummy{\noexpand\SLdummy}%
                \edef\SLoboxinhalt{}\edef\SLuboxinhalt{}%
                \SLobenfalse\SLuntenfalse
                \futurelet\next\SLsuchruf}
  \def\SLsuchruf{\ifx\next !\let\next\SLexpand
                 \else\let\next\SLerzeug\fi\next}
  \def\SLexpand#1#2#3{\ifx #2\sb\ifSLunten\let\SLspeicher\SLuboxinhalt
                   \else\def\SLspeicher{\SLtrick\SLabstand}\fi
                   \edef\SLuboxinhalt{%
                       \SLspeicher
                       \SLtrick\SLzwischen
                       \hbox\SLdummy{\hfil\mathsurround=0pt
$\SLtrick\scriptstyle\SLtrick#3$%
                                     \hfil}}%
                   \SLuntentrue%
                      \else\ifSLoben\let\SLspeicher\SLoboxinhalt
                   \else\def\SLspeicher{\SLtrick\SLabstand}\fi
                   \edef\SLoboxinhalt{%
                       \hbox\SLdummy{\hfil\mathsurround=0pt
$\SLtrick\scriptstyle\SLtrick#3$%
                                     \hfil}%
                       \SLtrick\SLzwischen
                       \SLspeicher}%
                   \SLobentrue\fi\futurelet\next\SLsuchruf}
  \def\SLerzeug{\def\SLtrick{}
                \setbox\SLhilfbox=\hbox{$\displaystyle{E}$}%
                \SLrandausgleich=.04\wd\SLhilfbox
                      \setbox\SLhilfbox=%
                         \hbox{\hskip -1\SLrandausgleich
                          \mathsurround=0pt$\displaystyle{\SLkern}$%
                               \hskip -1\SLrandausgleich}%
                      \SLhoehe = \ht\SLhilfbox
                      \advance\SLhoehe by \dp\SLhilfbox
                      \SLeffbreite = \wd\SLhilfbox
                      \advance\SLeffbreite by \SLab\SLhoehe
                      \ZBOX\SLubox{\vbox{\offinterlineskip
                                         \SLuboxinhalt
                                         \hrule height 0pt}}\SLdummy
                      \ZBOX\SLobox{\vbox{\offinterlineskip
                                         \SLoboxinhalt
                                         \hrule height 0pt}}\SLdummy
                      \SLsp = \SLzu\SLhoehe
                      \advance\SLsp by -.5\SLeffbreite
                      \SLuvorschub = -1\SLsp
                      \advance\SLuvorschub by -.5\wd\SLubox
                      \SLovorschub = -1\SLsp
                      \advance\SLovorschub by -.5\wd\SLobox
                      \advance\SLovorschub by .26\SLhoehe
                      \ifdim\SLuvorschub > \SLovorschub
                         \SLsp = \SLovorschub
                      \else
                         \SLsp = \SLuvorschub
                      \fi
                      \ifdim\SLsp < 0pt%
                         \advance\SLuvorschub by -1\SLsp
                         \SLmvorschub = -1\SLsp
                         \advance\SLovorschub by -1\SLsp
                      \else
                         \SLmvorschub = 0pt
                      \fi
                      \setbox\SLergebnis = \hbox{%
                         \offinterlineskip
                         \hskip\SLrandausgleich\relax
                         \vbox to 0pt{%
                            \vskip -1\ht\SLobox
                            \vskip -1\ht\SLhilfbox
\hbox{\hskip\SLovorschub\copy\SLobox\hfil}%
                            \hbox{\hskip\SLmvorschub\copy\SLhilfbox
                                  \hfil}%
\hbox{\hskip\SLuvorschub\copy\SLubox\hfil}%
                            \vss}%
                         \hskip\SLrandausgleich}%
                      \SLsp = \dp\SLhilfbox
                      \advance\SLsp by \ht\SLubox
                      \dp\SLergebnis = \SLsp
                      \SLsp = \ht\SLhilfbox
                      \advance\SLsp by \ht\SLobox
                      \ht\SLergebnis = \SLsp
                      \box\SLergebnis{}}
  \def\SLpruef{\ifx\next\SLsuchende\def\SLzu{0}\def\SLab{0}%
                  \def\next##1\SLsuchende{\relax}%
               \else\let\next\SLvergl
               \fi\next}
  \def\SLvergl#1_#2_#3:{\def\SLv{#1}%
                        \ifx\SLkern\SLv\def\SLzu{#2}\def\SLab{#3}%
\def\next##1\SLsuchende{\relax}%
                        \else\def\next{\futurelet\next\SLpruef}
                        \fi\next}
\def\PKT#1{#1!^\SLpunkt}

\newbox\minusbox
\def\minus{\mathchoice{\minusarb\displaystyle}%
                      {\minusarb\textstyle}%
                      {\minusarb\scriptstyle}%
                      {\minusarb\scriptscriptstyle}}
  \def\minusarb#1{\setbox\minusbox=\hbox{$#1x$}%
                  \hbox{\hskip .10\wd\minusbox
                        \vbox{\hrule width .6\wd\minusbox
                                     height .07\wd\minusbox
                              \vskip.53\ht\minusbox}%
                        \hskip .10\wd\minusbox}}

\def\Basis{\STAPEL e!_\SLstrich}
\def\C{\STAPEL C!_\SLstrich!_\SLstrich}
\def\X{\STAPEL X!_\SLstrich!_\SLstrich}
\def\Cinv{\STAPEL C!_\SLstrich!_\SLstrich^{-1}}
\def\Cg{\STAPEL C!_\SLstrich!_\SLstrich!^\SLstrich}
\def\CD{\STAPEL C!_\SLstrich!_\SLstrich!^\SLdreieck}
\def\XD{\STAPEL X!_\SLstrich!_\SLstrich!^\SLdreieck}
\def\ECKMT{\STAPEL M!_\SLstrich!_\SLstrich!^\SLdreieck^T}
\def\CgD{\STAPEL C!_\SLstrich!_\SLstrich!^\SLstrich!^\SLdreieck}
\def\BD{{\STAPEL C!_\SLstrich!_\SLstrich!^\SLdreieck}{^{\JIv}}}
\def\BDn#1{{\STAPEL C!_\SLstrich!_\SLstrich!^\SLdreieck}{^{\JIv}_{#1}}}
\def\Ttil{\STAPEL T!_\SLstrich!_\SLstrich!^\SLtilde}
\def\rhotil{\JI\STAPEL\varrho!^\SLtilde\!}
\def\rhotilzurho{{\JI\frac{\STAPEL\varrho!^\SLtilde\!}{\JI\!\varrho}}}
\def\rhozurhotil{{\JI\varrho\over\JI\STAPEL\varrho!^\SLtilde}}
\def\Skal#1{\STAPEL {#1}}
\def\Vek#1{\STAPEL {#1}!_\SLstrich}
\def\Ten2#1{\STAPEL {#1}!_\SLstrich!_\SLstrich}
\def\F{\STAPEL F!_\SLstrich!_\SLstrich}
\def\Fg{\STAPEL F!_\SLstrich!_\SLstrich!^\SLstrich}
\def\e{\STAPEL e!_\SLstrich!_\SLstrich}
\def\b{\STAPEL b!_\SLstrich!_\SLstrich}
\def\D{\STAPEL D!_\SLstrich!_\SLstrich}
\def\I{\STAPEL I!_\SLstrich!_\SLstrich}
\def\I{\STAPEL I!_\SLstrich!_\SLstrich}
\def\L{\STAPEL L!_\SLstrich!_\SLstrich}
\def\CnOldWed{\overset{\circ{}}{\Big(\STAPEL{C}!_\SLstrich!_\SLstrich_2\Big)}
                  \mskip1mu\mskip1mu
                  \mskip1mu {}^{\hat{\copy\JIGAMMz}}
                  \mskip1mu {}_{\hat{\copy\JIGAMMz}}}
\def\CnOldWedred{\STAPEL{C}!_\SLstrich!_\SLstrich!^\circ{}_2}

\def\CD{\C!^\SLdreieck}
\def\CDn#1{\C!^\SLdreieck_#1}
\def\X{\Ten2 X}
\def\XDn#1{\X!^\SLdreieck_#1}
\def\Cauchy{\STAPEL \sigma!_\SLstrich!_\SLstrich}
\def\Cn#1{\C_{#1}}
\def\Ln#1{\L_{#1}}
\def\Fn#1{\F_{#1}}
\def\qtil{\STAPEL q!_\SLstrich!^\SLtilde}
\newcommand{\IntGtil}[1]
   {\displaystyle{{\Big.^{\widetilde{\mathcal{G}}}}
      \hspace{-1.5ex}\int #1 \ {\rm d}\widetilde{V}}}
\newcommand{\IntGtilN}[2]
   {\displaystyle{{\Big.^{\widetilde{\mathcal{G}}_{#1}}}
      \hspace{-1.5ex}\int #2 \ {\rm d}\widetilde{V}_{#1}}}
\newcommand{\IntGhat}[1]
   {\displaystyle{{\Big.^{\widehat{\mathcal{G}}}}
      \hspace{-1.5ex}\int #1 \ {\rm d}\widehat{V}}}
\newcommand{\nadel}[1]
   {\big(\hspace{-3pt}\big(\hspace{1pt} #1\hspace{1pt}\big) \hspace{-3pt}\big)}
  
\newcommand{\inv}{{\minus 1}}
\newcommand{\invT}{{\minus\T}}    

\section{Introduction}
\label{sec:Intro}

For about 25 years, the clinical procedure of vertebroplasty is successfully applied to the treatment of osteoporotic vertebral bodies that have undergone a compression fracture. During this procedure, a biomaterial gets injected into the affected vertebral body and therein cures to a solid. This leads to an immediate reduction in pain as well as a stabilisation and strengthening of the spinal column \cite{HeiniEtal2001, JensenEtal1997}. 
Even though vertebroplasty is a well-accepted procedure, it may be accompanied by risks of complications caused by the treatment procedure itself or by the utilisation of specific injectable biomaterials. Regarding the treatment of vertebroplasty, a possible complication is cement leakage into unwanted regions of the vertebra. Furthermore, treated vertebral bodies exhibit an altered mechanical behaviour due to different mechanical properties of the biomaterial compared to human cancellous bone. This leads to different load distributions which may affect the treated vertebra itself but also adjacent vertebral bodies \cite{HeiniEtal2001,BaroudEtal2006}. 

There are different classes of injectable biomaterials used in vertebroplasty. However, acrylic bone cements play the most important role \cite{WebbSpencer2007}. This class of materials is a biocompatible but non resorbable polymer based on polymethylmethacrylate (PMMA). The application of acrylic bone cements within the scope of vertebroplasty can provoke further complications. Firstly, risk of thermal necrosis exists due to the exothermic chemical reaction which leads to a heat production and heating of the material. Moreover, release of remaining monomer caused by an incomplete chemical reaction may have a toxic impact on human tissue \cite{DahlEtal1994}.
 
Besides in vivo and in vitro measurement techniques, the process of vertebroplasty is subject of intensive research by numerical investigations. Generally, the overall procedure can be structured into three parts: the injection process during the operation, the curing process of bone cement inside the human body, and finally the long term behaviour of treated vertebrae. The injection process, which, from a physical point of view, is a fluid-dynamics problem, has been numerically examined using different computational methods. One basic approach is to investigate the injection process on the scale of a whole vertebral body. For this purpose, computational fluid dynamics (CFD) simulations are employed, as for example given by Teo \etal \cite{TeoEtal2007} or Widmer and Ferguson \cite{WidmerFerguson2011}. The simulations are based on the solution of Darcy's law in conjunction with region-specific intrinsic permeability, whereat specific values are directly connected to computer tomography data of the trabecular structure. In order to determine specific relations between trabecular bone morphology and corresponding intrinsic permeability, different approaches on the microstructural scale are reported. Baroud \etal \cite{BaroudEtal2004}, for instance, used a combination of experimental and analytical methods. Furthermore, numerical simulations of microstructural fluid flow through cancellous bone structure are utilized by using different numerical methods, like the lattice Boltzmann method  \cite{ZeiserEtal2008}, the finite volume method \cite{TeoTeoh2012} or smoothed-particle hydrodynamics \cite{BasafaEtal2010}. Additionally, morphological models are employed to analyse the microstructural flow of bone cement \cite{WidmerFerguson2013}.

The second stage of vertebroplasty, i.e. the curing of acrylic bone cement inside the human body, has been rarely investigated by numerical analyses. As mentioned above, the curing process of acrylic bone cements is accompanied by chemical shrinkage and an exothermic reaction, which leads to heating of the material and the vertebral body. The impact of bone cement shrinkage on the residual stresses of bone has, for example,  been numerically investigated by Kinzl \etal \cite{KinzlEtal2012}. In this study, the chemical shrinkage was adopted by a prescribed temperature evolution and a coefficient of thermal expansion. Thus, the cure kinetics have not been taken into account. Another issue of the curing process, i.e. the risk of thermal necrosis due to bone cement heating, has for example been analysed by Sta\'{n}czyk and van Rietbergen \cite{StanczykvanRietbergen2004}. Thereby, microstructural models have been employed for the simulation of temperature evolution during the curing process. However, the mechanics are not considered in this study. Within the works mentioned in this paragraph, microstructural finite element (FE) models are employed. A common approach of obtaining such models is given by the $\mu$FE technologies \cite{vanRietbergenEtal1995}. 

Among numerical studies currently available in the literature, the major focus is on the simulation of the long term behaviour of treated vertebrae by FE-simulations. Thereby, several studies simply assume an approximate bone cement distribution (see, for example, \cite{RohlmannEtal2010} for the case of vertebroplasty and \cite{LiebschnerEtal2001} for the case of kyphoplasty) or employ a direct conversion of augmented vertebral bodies to FE-models by imaging methods \cite{ChevalierEtal2008}. An overview over different models is for example given in \cite{JonesWilcox2008} or \cite{JonesEtal2012}. A different approach is to combine injection simulations and FE-simulations of vertebral body loading, as for example given by the above mentioned study \cite{TeoEtal2007}. However, this study does not take into account the curing process itself and the temperature evolution during bone cement curing. 

The aim of this paper is to present a new approach of numerical simulation which combines several of the above mentioned simulation strategies and, moreover, includes an enhanced phenomenological description of the material behaviour of acrylic bone cements. Thus, a combination of thermo-mechanically coupled simulations of the different stages of vertebroplasty is presented. More precisely, this integrated modelling and simulation framework includes 

\begin{itemize}
   \item[(i)] the generation of microstructural computer models based on $\mu$CT images of human cancellous bone,
   \item[(ii)] the CFD-simulation of acrylic bone cement injection into the cancellous bone structure, and
   \item[(iii)]  the FE-simulation of bone cement curing inside the trabecular structure.
\end{itemize}

With the help of this approach, detailed investigations on the microstructure of trabecular bone in conjunction with acrylic bone cements are facilitated. This includes studies on the bone cement's flow behaviour during injection, the impact of the exothermic chemical reaction and the related heat dissipation on human tissue, as well as arising stresses due to volume changes of the bone cement. Moreover, an analysis of stiffness properties of treated and untreated cancellous bone is enabled.  An essential part of our approach is the proper representation of the material behaviour of acrylic bone cements during the different simulation stages. Thus, an outline of the addressed phenomenology and corresponding material modelling approaches will be given.

The paper is organised as follows. In Section \ref{sec:bone_cement}, an overview of the phenomenology and corresponding modelling approaches for acrylic bone cements is provided. Next, the integrated modelling and simulation framework is introduced in Section \ref{sec:framework}. Finally, Section \ref{sec:testcase} exemplifies the application of the framework to a sample of human cancellous bone. Thereby, the injection of initially liquid bone cement and the in vivo curing of the material are considered. Besides the thermo-chemical behaviour during cure, the evolving stresses in cancellous bone due to chemical shrinkage are investigated and the influence of boundary conditions and bonding behaviour at the interface between bone and bone cement are addressed.

\section{Acrylic bone cements: Phenomenology and material modelling approaches}
   \label{sec:bone_cement}
\subsection{Phenomenology of acrylic bone cements}
  \label{sec:bone_cement_phaeno}
	
\begin{figure}[ht]
  \centering
    \includegraphics[width=0.4\textwidth]{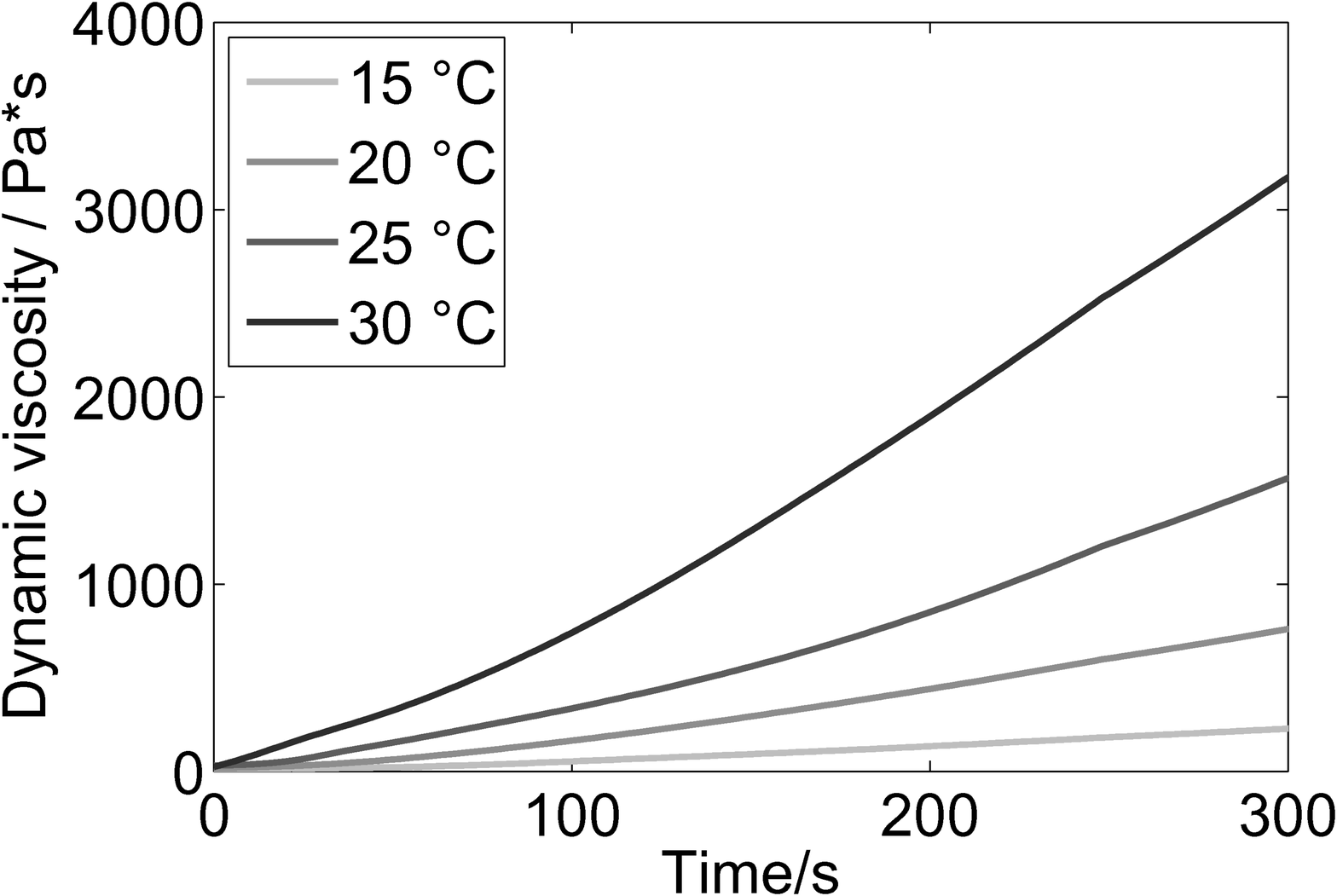}
	\hspace{8mm}
    \includegraphics[width=0.4\textwidth]{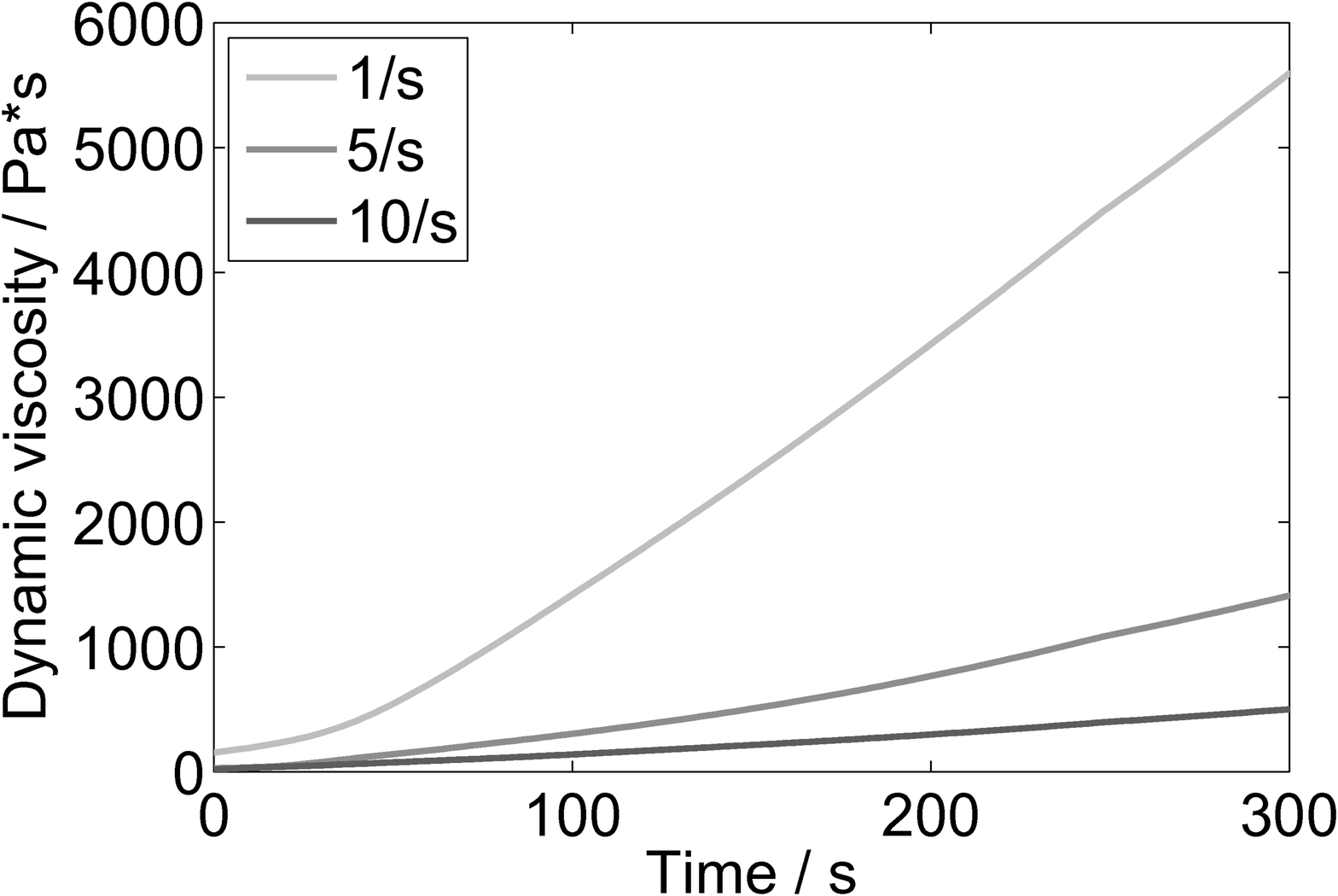}
    \caption{Viscosity at constant shear rate $\dot{\gamma} = 5 \,s^{-1}$ (left) and at constant temperature $\theta = 25^{\circ}C$ (right).}
    \label{fig1}
\end{figure}

Acrylic bone cements are based on the thermoplastic polymer polymethylmethacrylate (PMMA). In clinical applications, the material is mixed out of two components, a polymer powder and a liquid monomer, whereat the components include different additives, e.g. radio-pacifiers among others \cite{WebbSpencer2007}. Immediately after mixing two different processes start. Firstly, a physical dissolution process (also referred to as swelling) occurs due to penetration of liquid monomer into the polymer powder. This leads the two components to form a liquid or rather a dough. In this state, acrylic bone cements show distinctive non-Newtonian flow behaviour \cite{FarrarRose2001, KrauseEtal1982}. This fact has been investigated by detailed rheological measurements on a cone and plate rheometer. As given in Fig. \ref{fig1}, three major influences on the viscosity can be observed. Firstly, the viscosity rises with elapsing time, meaning with proceeding dissolution. Furthermore, an increase in temperature leads to higher values for dynamic viscosity (see Fig. \ref{fig1}a). In fact, higher temperatures accelerate the progress of dissolution. This effect of dissolution on the viscosity gets superposed by a shear thinning behaviour of acrylic bone cement material, i.e. increasing shear rate causes a decrease in viscosity (see Fig. \ref{fig1}b).

\begin{figure}[ht]
 \includegraphics[width = 0.55\textwidth]{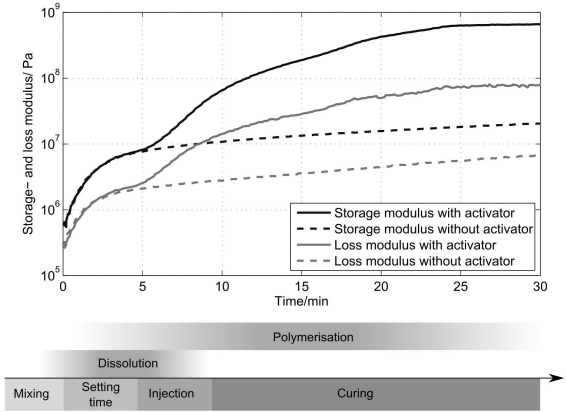}
 \caption{Storage and loss modulus of a standard bone cement and a bone cement composition without activator (no polymerisation). The different curves reveal the different impact of the processes of physical dissolution and polymerisation.}
 \label{fig2}
\end{figure}

Besides the dissolution process, the chemical process of radical polymerisation starts immediately after mixing as well. Thereby, activator and initiator, which are also included in the powder and the liquid, form free radicals and lead to a lengthening of the polymer chains \cite{KuehnEtal2005}. The polymerisation process is responsible for the curing of the material, i.e. the initially liquid material undergoes a phase change to a solid.  Even though both dissolution and polymerisation are simultaneous, it can be observed that the dissolution process has a more dominant effect at the beginning while the polymerisation process reaches more significance after a certain time. Fig. \ref{fig2} depicts this behaviour in the case of mechanical testing on a parallel-plate rheometer. Thereby, a commercial acrylic bone cement including initiator and activator is compared to a bone cement mixture where the polymerisation is not initiated. Thus, the different impacts of dissolution and polymerisation on the mechanical stiffness can be observed. Furthermore, Fig. \ref{fig2} indicates the viscoelastic behaviour in the solid state of the material. 

The polymerisation process has further consequences. Firstly, it is exothermic and therefore is accompanied by heat production. This typically leads to temporary temperature peaks under natural environmental conditions, which also occurs in application \cite{BelkoffMolloy2003}. Moreover, volume shrinkage of up to $6-7 \%$ can be observed for commercial acrylic bone cements \cite{KuehnEtal2005}, and further properties like the specific heat capacity and the coefficient of heat conduction change during the polymerisation process \cite{KolmederEtal2011}. 

In order to simulate the material behaviour of acrylic bone cements during the curing process, a detailed experimental characterisation as well as a phenomenological modelling approach including the coupling of chemical, thermal and mechanical processes is needed. In previous studies, the commercial acrylic bone cement Osteopal\textsuperscript{\textregistered}V manufactured by Heraeus Medical, Germany has been examined and characterised relating to its fluid behaviour \cite{KolmederLion2013} and during the curing process \cite{KolmederLion2010, KolmederEtal2011}. The mechanical behaviour of fully cured bone cement and a corresponding viscoelastic model at finite strains has been presented Lion \etal \cite{LionEtal2008}.  In the following, concrete material modelling approaches for the representation of acrylic bone cements as a non-linear fluid and as a non-linear curing dependent viscoelastic solid are provided (see subsequent sections \ref{sec:fluid} and \ref{sec:solid}, respectively). 

\subsection{Fluid flow model}
  \label{sec:fluid}

The considered acrylic bone cements exhibits non-linear fluid flow behaviour \cite{KolmederLion2013}. In order to represent this behaviour within a phenomenological model, firstly an internal variable $p$, representing the progress of the dissolution process, is introduced.  A first order differential equation is chosen for the evolution of this internal variable:
\begin{equation}\label{eq:disso}
  \dfrac{{\rm d}\,p}{{\rm d}t}
	=  \begin{cases}
            \left(k_{p1} + k_{p2}\, p^m\right),  & \text{if }p<p_{max}    \\
            0,                                  & \text{if }p\geq p_{max}
      \end{cases} \ .
\end{equation}
The rate of $p$ is distinguished by cases to stop the evolution, once a maximum value for $p$ is reached. For $p<p_{max}$, the coefficients $k_{p1}$ and $k_{p2}$, showing an Arrhenius behaviour, take into account the temperature dependence of the dissolution process as follows
\begin{equation}\label{eq:arrhenius}
     k_{pi} = 
		    A_{pi}\cdot \exp\left(\dfrac{-E_{pi}}{R\,\theta}\right) 
		\ \text{;} \ i = 1,2 \ .
\end{equation}
Here, $R = 8.314 \, \rm J/K mol$ is the universal gas constant, $\theta$ is the thermodynamic temperature and $E_{pi}$ are activation energies. Furthermore, $A_{pi}$ are constant material parameters. To capture the shear thinning behaviour, a power law fluid model is applied. In combination with a function that depends on the internal variable $p$, the newly developed model \eqref{eq:viscosity} for the viscosity $\eta$ is able to cover all previously described effects. 
\begin{equation}\label{eq:viscosity}
 \eta = \eta_0(p) \, l \left|\dfrac{\dot{\gamma}}{\dot{\gamma_0}}\right|^{n-1} \ , \quad \eta_0(p) = h_{\eta}p
\end{equation}
Therein, $\dot{\gamma}$ is the shear rate. In a first approach, the function $\eta_0(p)$ is set to $h_{\eta}p$, whereat $h_{\eta}$ is a proportionality constant. To preserve units, $\dot{\gamma_0}$ with a value of $1 \,s^{-1}$ is introduced. The identified parameters of Eqs. \eqref{eq:disso}, \eqref{eq:arrhenius} and \eqref{eq:viscosity} are listed in Table~\ref{tab1}.

As mentioned in section \ref{sec:bone_cement_phaeno}, polymerisation starts right away after mixing. Analogous to the progress of dissolution, another internal variable $q$ (referred to as degree of cure) is introduced for the representation of the progress of polymerisation. The evolution is also described by a first order differential equation:
\begin{equation}\label{eq:polymerisation}
  \dfrac{{\rm d}\,q}{{\rm d}t}= (k_{q1} + k_{q2} \,  q^\alpha)\,(1-q)^\beta \, f_D(q,q_{end}) \ .
\end{equation}
The values of the degree of cure are defined in the range of $q \in [0,1]$, whereat $q = 1$ represents the fully polymerised material. Again, $k_{q1}$ and $k_{q2}$ show Arrhenius behaviour. It has to be noted, that the evolution strongly depends on temperature and a value of $q = 1$ will not be reached if the temperature is below the glass transition temperature. This behaviour is governed by the function $f_D(q,q_{end})$ which in turn depends on a function $q_{end}$ representing the maximum of the degree of cure to be reached in isothermal conditions. More details on the experimental characterisation and the adjustment of model parameters in Eq. \eqref{eq:polymerisation} to the acrylic bone cement Osteopal\textsuperscript{\textregistered}V can be found in \cite{KolmederLion2010}.

\begin{table}[ht]
  \caption{Identified material parameters.} 
  {\begin{tabular}{p{1.5cm}p{2.5cm}p{1.5cm}p{3cm}p{1.5cm}p{1.3cm}}
    \hline
    Parameter & Value & Parameter & Value & Parameter & Value\\ 
    \hline \\[-3mm]
      $\ \ A_{p1}$          & $3.0 \times 10^{14}\,s^{-1}$                 
	    & $\ \ E_ {p1}$    & $8.232 \times 10^{4}\;\text{J}/\text{mol}$    
        & $\ \ h_{\eta}$   & $1\,\text{Pa}\,\text{s}$   \\
	  $\ \ A_{p2}$          & $0.0023\,s^{-1}$              
        & $\ \ E_{p2}$     & $312.6\;\text{J}/\text{mol}$   
		& $\ \ l$              & $8.996$   \\  
      $\ \ m$                 & $1$                                          
	    & $\ \ p_{max}$   & $2000$           
	    & $\ \ n$             & $0.2601$                       \\ 
    \hline
  \end{tabular}}
  \label{tab1}
\end{table} 

\subsection{Viscoelastic solid including polymerisation}
  \label{sec:solid}

For the consideration of bone cement curing from a dough to a solid, a modelling approach within the framework of non-linear continuum mechanics is employed. Thereby, the following aspects as well as their couplings are taken into account:

\begin{itemize}
   \item the polymerisation process (given by Eq. \eqref{eq:polymerisation}),
   \item the process-depending viscoelastic behaviour during polymerisation (cf. Fig. \ref{fig2}),
   \item the viscoelastic behaviour of the polymerised material,
   \item the transient heat transfer behaviour, and
   \item volume changes (heat expansion and chemical shrinkage).
\end{itemize}

The kinematics are described by the deformation gradient $\F$ which gets multiplicatively decomposed into a thermo-chemical part $\F_{\theta C}$ and a mechanical part $\F_M$ as follows:
\begin{equation}\label{eq:FThetaC}
  \F = \F_M\cdot\F_{\theta C}
  \qquad \quad
  \F_{\theta C} = \sqrt[3]{J_{\theta C}(\theta, q)} \, \underline{\underline{I}} \, .
\end{equation}
The thermo-chemical volume changes are assumed to be isotropic. Thus, $\F_{\theta C}$ is described by the scalar volume ratio  $J_{\theta C}(\theta,q)$ which is a function of the temperature $\theta$ and the degree of cure $q$. More precisely, the following ansatz has been chosen:
\begin{equation}\label{eq:JthetaC}
  J_{\theta C}(\theta, q) = {\rm exp}\Big[\big(\alpha + \Delta \alpha \, q\big)\big(\theta - \theta_{ref}\big) + \Delta \beta \, q\Big] \ .
\end{equation}
The material parameters $\alpha, \Delta \alpha$ and $\Delta \beta$ have been identified using an experimental setup based on the Archimedes' principle \cite{KolmederLion2010}.

In order to obtain a nearly incompressible formulation of the mechanical part of the material model, a further decomposition is employed. Thus, the mechanical part of the deformation is decomposed into a volumetric deformation $J_M = {\rm det} \, \F_M$ and a volume preserving part $\F_G$. Furthermore, a symmetric stretch tensor $\C_G$ of right Cauchy-Green type is defined and will be employed as input for the subsequent formulation. 
\begin{equation}\label{eq:FG}
    \F_G = \dfrac{1}{J_M^{1/3}} \ \F_M
    \qquad 
    \C_G = \F_G^T\cdot \F_G
\end{equation}

Next, a  strain-energy function (Helmholtz free energy density) is formulated, which describes the energy storage of the material. It is additively decomposed into a pure mechanical volumetric part $\bar{\psi}_{vol}$, an isochoric viscoelastic part $\bar{\psi}_{vis}$ and a thermo-chemical part $\bar{\psi}_{\theta C}$. 
\begin{equation}\label{eq:FreeEnergy}
  \bar{\psi} = \bar{\psi}_{vol}(J_M) + \bar{\psi}_{visc}(\C_G,z) + \bar{\psi}_{\theta C}(\theta, q)
\end{equation}
Thereby, the pure mechanical volumetric part $\bar{\psi}_{vol}(J_M)$ ensures the compressibility of the material and solely depends on the mechanical volume ratio $J_M$ and the bulk modulus $K$. In this work, the ansatz
\begin{equation}\label{eq:psiVol}
  \bar{\psi}_{vol} = \dfrac{9\,K}{2} \big( J_M^{1/3} -1\big)^2
\end{equation}
has been employed. However, other relations may be used.
The second part of Eq. \eqref{eq:FreeEnergy} contains a formulation of a chain of several Maxwell-elements connected in parallel:
\begin{equation}\label{eq:psiVisc}
   \bar{\psi}_{visc} = \sum_{n=1}^{N}  \bar{\psi}_{visc,n}.
\end{equation}
The strain-energy of a single Maxwell-element reads as
\begin{equation}\label{eq:psiViscK}
 \bar{\psi}_{visc,n} 
     = \Bigg\{- \displaystyle\int\limits_{-\infty}^{z} G_n(z-s)  \ \bigg( \dfrac{\rm d}{{\rm d}s}\C_G^{-1}(s)\bigg){\rm d}s \Bigg\} \cdot \hspace{-2pt} \cdot \hspace{2pt} \C_G 
\end{equation}
with the kernel function
\begin{equation}\label{eq:psiViscKKernel}
   G_n(z-s) =  2 \, \mu_n {\rm exp}\Big[-\dfrac{z-s}{\tau_n}\Big]  \ .
\end{equation}
Therein, $\mu_n$ and $\tau_n$ denote the stiffness parameter and the relaxation time, respectively,  for the $n^{th}$ Maxwell-element \cite{HauptLion2002}. Moreover, the formulation is extended by a variable $z$ which enables the consideration of process depending properties of the Maxwell-elements. More precisely, this variable is formulated as an ordinary differential equation that depends on the temperature $\theta$ and the degree of cure $q$:
\begin{equation}\label{eq:intrinsic}
  \dfrac{{\rm d}\, z}{{\rm d}t} = \dot{z} = f(\theta,q) \ .
\end{equation}
In this manner, the temperature and degree of cure dependencies of the viscoelastic properties are taken into account. A specific ansatz can be found in \cite{KolmederEtal2011}. The remaining term $\bar{\psi}_{\theta C}(\theta, q)$ of Eq. \eqref{eq:FreeEnergy} represents the thermo-chemical stored energy which has been formulated and identified based on differential scanning calorimetry (DSC) measurements \cite{LionYagimli2008,KolmederEtal2011}.

Following the common approaches in non-linear continuum mechanics, the next steps are the evaluation of the free energy within the Clausius-Duhem inequality to obtain a thermodynamically consistent material behaviour. Furthermore, the $1^{\rm st}$ law of thermodynamics is employed to account for the thermo-mechanical coupling. This procedure results in a system of non-linear equations describing the fully thermal-chemical-mechanically coupled material behaviour. Among others, it contains Eq. \eqref{eq:stresses} for the calculation of stresses of second Piola-Kirchhoff type as well as evolution equations for the degree of cure $q$  and the material function $z$ (see Eqs. \eqref{eq:polymerisation} and \eqref{eq:intrinsic}, respectively). 
\begin{equation}\label{eq:stresses}
\begin{array}{c}
   \Ttil = \Big[ \, \Ttil_G \cdot \C_G\Big]' \cdot\C^{-1} 
           \, + \, 3 \, K \, J_M^{1/3} \, (J_M^{1/3}-1) \; \C^{-1} \\[4mm]
   \displaystyle
   \Ttil_G = - \sum_{n=1}^{N} 
               \ \int\limits_{-\infty}^{z} G_n(z-s) \  
                   \bigg( \dfrac{\rm d}{{\rm d}s}\C_G^{-1}(s)\bigg)
               {\rm d}s
\end{array}
\end{equation}
Here, the superscript $\big(...\big)'$ represents the deviatoric part of a tensor of rank 2. The second part of the system of equations is given by an equation of transient heat conduction:
\begin{equation}\label{eq:heatcond}
  c_{\theta} \, \dot{\theta} - \widetilde{\underline{\nabla}} \cdot \qtil + \mathcal{D}_{\theta C} + \mathcal{D}_{vis} - w_{e} = 0 \ .
\end{equation}
Therein $c_{\theta}$ is the heat capacity per unit volume. The rates of dissipated energy due to exothermic chemical reaction of polymerisation and the viscoelastic material behaviour are denoted by $\mathcal{D}_{\theta C}$ and $\mathcal{D}_{visc}$ , respectively. The function $w_e$ represents the thermo-elastic coupling. All of those mentioned functions are derived by the free energy function \eqref{eq:FreeEnergy} in conjunction with the first law of thermodynamics. Additionally, $\qtil$ denotes the heat flux vector on the reference configuration for which the Fourier's law has been employed.

The outlined material model has been implemented into the FE-software MSC.MARC\textsuperscript{\textregistered} and material parameters have been fitted to the bone cement Osteopal\textsuperscript{\textregistered}V. To this end, different methods of experimental characterisation and automated, non-linear parameter identification have been employed \cite{KolmederLion2010,KolmederEtal2011}. An application of the material model within FE-simulations has been presented in \cite{LandgrafEtal2012}.

\section{Integrated framework of model generation and numerical simulation}
  \label{sec:framework}

\begin{figure}
  \centering
 \includegraphics[width = .99\textwidth]{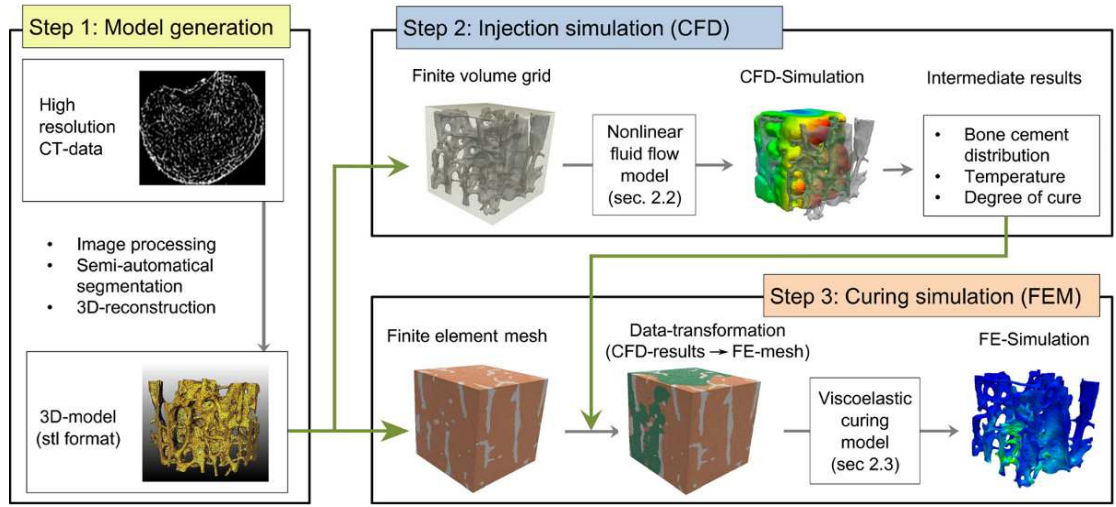}
 \caption{Schematic representation of the modelling and simulation framework for the computation of bone cement injection and in vivo curing.}
 \label{fig3}
\end{figure}

As outlined in Section \ref{sec:Intro}, the aim of the present work is the simulation of thermo-mechanical processes related the surgical technique of vertebroplasty. This includes the injection of bone cement into vertebral cancellous bone as well as the in vivo curing of the injected biomaterial. Moreover, studies on  the stiffness properties of treated and
untreated cancellous bone are enabled. For this, an integrated modelling and simulation framework was built up (cf. Fig. \ref{fig3}) which consists of the following three main steps: 
\begin{enumerate}
 \item Generation of the 3D-surface-model of an osteoporotic trabecular structure based on $\mu$CT-data,
 \item Injection simulation by CFD-methods, and
 \item FE-simulation of curing and post-vertebroplasty processes.
\end{enumerate}

Firstly, the generation of computer models of osteoporotic cancellous bone structures is carried out using greyscale-data in DICOM-format obtained by $\mu$CT-scanning technology. This is performed by means of an iterative procedure of image processing, segmentation and 3D-modelling techniques (see Section \ref{sec:ImageProc}). The final computer models of osteoporotic cancellous bone are provided in stereolithography file format (STL) and are utilised as input for the subsequently described CFD- and FE-models.

The second step incorporates  the CFD-simulation of bone cement injection into the trabecular structure of a vertebral body. Thereby, the just mentioned STL-model is employed as input. In order to simulate the non-linear and process-depending material behaviour of acrylic bone cements during the injection process, the fluid flow model described in  Section \ref{sec:fluid} is employed. The simulation is carried out using a thermo-mechanically coupled solver for a two-phase-flow of immiscible fluids. Thus, the flow behaviour of bone marrow in considered as well. More information are provided in Section \ref{sec:CFD}.

Subsequently to the injection process, the curing of the acrylic bone cement inside the vertebral body is considered by non-linear FE-simulations (see Section \ref{sec:FEM_modelling}). Again, the STL-model of cancellous bone is used as input for the development of the FE-mesh of both bone tissue and the pores in between. Furthermore, the bone cement distribution and state variables like the temperature and the degree of cure are transferred from the CFD-simulation to the FE model. The thermo-mechanically coupled curing processes of bone cement are described by the viscoelastic material model presented in Section \ref{sec:solid}.  Thus, the effects of curing, volume changes and evolving stiffness of the material can be analysed.

In the following, the three main steps of the framework are described in more detail.

\subsection{Generation of 3D-surface models of osteoporotic trabecular structures}
  \label{sec:ImageProc}

In our study, the generation of 3D-surface models of osteoporotic cancellous bone is based on greyscale-data in DICOM-format, obtained by $\mu$CT-scanning technology. One famous approach for the model generation process is given by so-called $\mu$FE technologies \cite{vanRietbergenEtal1995}.  These techniques have the advantage of delivering models of approved accuracy (e.g. \cite{UlrichEtal1998}), but in turn are of very high complexity and typically consist of a large number of elements. As the study described in this article is dedicated to both highly non-linear CFD- and FE-simulations, we refer to an alternative approach similar to the one developed in \cite{KoberEtal2007}. Notably, Ulrich \etal \cite{UlrichEtal1998} already showed that model generation of acceptable accuracy with reduced element number is possible. 

For the subsequent mesh generation for CFD- and FE-simulations, different requirements on the computer models have to be fulfilled. Firstly, the 3D-surface models have to be fully intersection-free and manifold.\footnote{From a mathematical point of view, a manifold surface is a topological space, where every point has a neighbourhood topologically equivalent to an open disc or half disc of $R^2$, see inter alia \cite{YuEtal2010}.} Thus, the surface must not contain any connections by ``free edges'' or ``free points''.  Fulfilling standard quality criteria as smallest dihedral angle above $5$ (better $10$) degrees, the 3D-surface model's complexity should not exceed $1000$ elements per $\rm mm^3$. Furthermore, it should be connected, so comprising only one component. This, for example, is not evident for a sample cut out of some vertebral body. In order to fulfil the given requirements, the procedure described below is applied. Thereby, all image processing, segmentation, and 3D-reconstruction steps are performed using the visualisation and modelling tools available in the software Amira\textsuperscript{\textregistered}.\footnote{VSG - Visualization Sciences Group, Berlin, Germany,  http://www.amira.com/}

Firstly, the DICOM images are loaded to the software and visualised in multiplanar viewing (axial, coronal, and sagittal) for the sake of orientation and selection of an appropriate region of interest. As kind of preprocessing, repeated Gaussian smoothing of the $\mu$CT-images and analysis of the grey values along the single trabeculae are performed. By 3D-region growing, the trabecular tissue is semi-automatically segmented from the $\mu$CT-images. After moderate Gaussian smoothing of the segmentation, a connected component is "picked" from the segmentation in order to get a connected surface. 

Immediate 3D-reconstruction with some marching cube like algorithm generally delivers models by far exceeding the complexity acceptable for the present study. Furthermore, subsequent simplification, including necessary smoothing, produces plenty of intersections and non-manifold contacts. This is cured by means of an iterative procedure similar to the one introduced in \cite{KoberEtal2007}. Within one iteration loop, the model is reconstructed in 3D, then stepwisely simplified and smoothed, until it has the appropriate number of elements with, at the same time, acceptable smoothness and appearance. If necessary and feasible, the model is also interactively edited. By intersection of the surface model with the original grey images, it is assured that the model is still in agreement with the CT-data. Thereafter, the modified surface model is converted again into some labelled image stack with the same dimensions as the original data. By this procedure, the regularity of the new labelled image stack is iteratively increased and, as a result, the number of intersections and non-manifold contacts of the respective 3D-surface model are reduced. The iteration is stopped once the number of intersections and non-manifold contacts is so small that it can be interactively removed. Finally, semi-automatic improvement of triangular quality completes the modelling process and closed 3D-models are exported to STL-file format, fully manifold and of suitable element quality. Those models are employed for the definition of the CFD-grid and the mesh of the FE-model.

\subsection{Finite volume modelling and implementation}
\label{sec:CFD}

The open source computational fluid dynamics toolbox OpenFOAM\textsuperscript{\textregistered}\footnote{OpenCFD Ltd (ESI Group), Wokingham, UK, http://www.openfoam.com/} is used to simulate the injection process. Advantages of the software are full accessibility of the source code and the possibility of parallelisation. Moreover, it provides a mesh generation utility snappyHexMesh that is able to generate meshes automatically from triangulated surface geometries in STL-file format. For this purpose, an initial hexahedral mesh is needed which defines the geometrical size of the computational domain and a base level mesh density. The snappyHexMesh utility provided by OpenFOAM\textsuperscript{\textregistered} refines the background mesh around features and surfaces specified in the provided STL-geometry, removes cells within unwanted regions, e.g. inside the bone tissue, and creates a castellated mesh. A further process step snaps the cells to the specified surface. 

To account for the specific fluid flow behaviour of acrylic bone cements within the CFD-simulation, the mathematical fluid flow model described in section \ref{sec:fluid} has been included into the software. Here, it has to be noted that the vertebral body is filled with bone marrow prior to the injection of bone cement. Thus, two different fluids have to be considered within the simulation of the injection process. Towards this end, the implementation of the fluid flow model is carried out on the basis of the incompressible volume-of-fluid solver interFoam \cite{Rusche2002}, a solver for two immiscible fluids.  The solver offers a phase-fraction based interface capturing method to calculate the boundary surface between bone cement and bone marrow \cite{NohWoodward1976, HirtNichols1981}. Thereby, the phase fraction $\alpha_p$ is introduced as follows:
\begin{equation}\label{eq:phasefraction}
  \alpha_p \in [0,1]
    :\ \left\{\begin{array}{cl}
        \alpha_p = 1      & \dots \text{bone cement cells} \\      
        0<\alpha_p<1 & \dots \text{interfacial cells} \\        
         \alpha_p =0      & \dots \text{bone marrow cells}         
      \end{array} \right.
\end{equation}
The variable $\alpha_p$ defines the volume fraction of the bone cement phase within a finite volume cell. Thus, $\alpha_p=1$ describes a finite volume cell containing only bone cement. A value of $\alpha_p = 0$ denotes a cell containing only bone marrow.

To take into account the special behaviour of bone cement, several modifications have been made to the standard solver interFoam. In general, three different physical phenomena can be identified that interact with each other and govern the behaviour of bone cement during the injection phase: Non-Newtonian fluid behaviour according to Eq. \eqref{eq:viscosity}, the dissolution process (see Eq. \eqref{eq:disso}) and transient heat transfer. The latter is added to the standard solver by implementing a further transport equation as follows:
\begin{equation}\label{eq:heattransfer}
  \dfrac{{\rm d}\rho h}{{\rm d}t} + \Vek{\nabla} \cdot \big(\rho h \, \Vek{v} -  k \, \Vek{\nabla}\,\theta\big) = \Ten2{\tau} \cdot\hspace{-.5ex}\cdot\hspace{.5ex} \D
  \ , \quad 
  \D = {\rm sym} \big(\Vek{\nabla} \otimes \Vek{v}\big) \ .
\end{equation}
Therein, $h$ is the specific enthalpy, $\rho$ the density, $k$ the thermal conductivity, $\theta$ the thermodynamic temperature, $\Ten2{\tau}$ the deviatoric part of the Cauchy stress tensor and $\Vek{v}$ the fluid velocity. The term on the right hand side of Eq. \eqref{eq:heattransfer}$_1$ represents the stress power which can also be interpreted as viscous dissipation.

In order to keep track of the dissolution and polymerisation processes throughout the solution domain, Eqs. \eqref{eq:disso} and \eqref{eq:polymerisation} have to be extended. Thereby, dissolution and polymerisation are supposed to be transported only by convection, but not by diffusion, and are only valid in the bone cement phase. In the case of the dissolution process, described by the evolution equation \eqref{eq:disso}, the extension by a convection term reads as
\begin{equation}\label{eq:disso_conv}
  \dfrac{{\rm d}p}{{\rm d}t} + \Vek{\nabla} \cdot (p \, \Vek{v})= (k_{p1} + k_{p2} p^m) \ .
\end{equation}
Eq. \eqref{eq:polymerisation} is adjusted in the same manner. Having the heat transfer equation \eqref{eq:heattransfer} and the progress of dissolution \eqref{eq:disso_conv} available in the program code, the viscosity model according to Eq. \eqref{eq:viscosity} can finally be added to the solution procedure.

\subsection{Finite element modelling and data transformation}
\label{sec:FEM_modelling}

Bone cement curing is investigated by non-linear, thermo-mechanically coupled FE-simulations. All modelling and simulation steps are accomplished using the preprocessor MSC.MENTAT\textsuperscript{\textregistered} and the non-linear FE-solver MSC.MARC\textsuperscript{\textregistered}.\footnote{MSC.Software Corporation, Santa Ana, California, USA, http://www.mscsoftware.com/}

In order to obtain a smooth transition from the CFD-simulation to the FE-model, finite element meshes are build up using the same input models as already used for the definition of the corresponding CFD-grids. Thereby, both bone tissue and pores are meshed by tetrahedrons.  Whilst the meshing of bone tissue can be accomplished fully automatically, some preprocessing steps have to be fulfilled for the meshing of the pores. This includes the definition of a region to be meshed (in our case this region is a cube) and the construction of a triangular mesh on the surfaces of the cube. By this, the holes of the pore volume on the boundary are closed. Subsequently, the closed volume can be meshed automatically.

Next, the initial state of the FE-simulation has to be defined. This state shall represent the same physical state as given by the CFD-simulation. Thus, the bone cement distribution and different variables like the temperature field and the internal variables resulting from the CFD-simulation have to be passed to the FE-mesh. However, since CFD-grid and FE-mesh generally do not coincide, a data transformation of the fluid simulation results to the FE-model is required. Towards this end, a converter has been developed that reads the result files provided by the fluid simulation software, converts the CFD-results to the FE-mesh, and writes appropriate input files for the FE-model. 

In the first step of conversion, the distinction between bone cement and bone marrow within the FE-mesh is accomplished. Here, a nearest neighbour search including a space partitioning scheme is fulfilled within a loop over every integration point of the finite elements (elements representing bone tissue are excluded). The algorithm detects the cell of the CFD-grid which is located nearest to the current integration point of the FE-mesh. Having found the nearest cell, the phase fraction of the cell (cf. Eq. \eqref{eq:phasefraction}) governs the material of the integration point, i.e. bone marrow if $\alpha_p < 0.5$ or bone cement if $\alpha_p \ge 0.5$.  Simultaneously, the temperature and, in the case of bone cement, the degree of cure are transferred from the chosen cell of the CFD-grid to the current integration point. In a further step, the nearest neighbour search algorithm is adopted to assign temperatures to all the nodes of the FE-mesh.

A crucial aspect of numerical implementation shall be noted at this point. The initial values for the temperature and the degree of cure, both resulting from CFD-analyses, generally differ from predefined values for the reference temperature $\theta_{ref}$ and the reference degree of cure $q_{ref} = 0$ (see Eq. \eqref{eq:JthetaC}). In such situations, the volume would tend to change its volume immediately. This may lead to unrealistic behaviour and, moreover, to mesh distortions right at the beginning of finite element analyses. In order to cope with this effect, a special numerical algorithm has been developed. It takes into account deviations of initial values from their reference values and corrects the volume ratio \eqref{eq:JthetaC} in such a manner that the initial volume of the FE-mesh remains undeformed (see \cite{LandgrafEtal2013}).

The finite element simulation of bone cement curing is accomplished using the staggered solution scheme provided by MSC.MARC\textsuperscript{\textregistered}. Thereby, the mechanical part of the solution is treated as a static problem, thus inertia effects are not taken into account, and the thermodynamic part is calculated by a transient solver. Furthermore, MSC.MARC\textsuperscript{\textregistered} provides different user interfaces to include self-written code.  The main user interface utilised in this work is the user subroutine \textit{HYPELA2} for the implementation of  the stress-strain relation given by Eq. \eqref{eq:stresses}. Additionally, the user subroutines \textit{UPSCHT} and \textit{FLUX} are employed for the definition of a user defined specific heat capacity $c_\theta$ and the rates of dissipated energy $\mathcal{D}_{\theta C}$, $\mathcal{D}_{visc}$ and $w_e$ given in Eq. \eqref{eq:heatcond}, respectively. Beyond that, the user subroutines \textit{IMPD} and \textit{INITSV} are adopted for the definition of the initial conditions of the FE-model.

\section{Application of the framework: A case study}
  \label{sec:testcase}

In this section, the integrated modelling and simulation framework, described in the preceding section, is applied to a case study consisting of a sample of trabecular bone structure.  To this end, high resolution axial $\mu$CT-data\footnote{vivaCT 40, Scanco Medical AG, Br\"uttisellen, Switzerland, data acquired at Department of Trauma Surgery/Biomechanics, Medical University Innsbruck, Austria} of a post-mortem preparation (female, 69 Y) were acquired with isotropic voxel size of $19~\mu$m in $x-$, $y-$, and $z-$direction and $1804 \times 1824$ image matrix. These data represent an osteoporotic thoracic (th10) vertebral body with the vertebral processes removed. The data set is in DICOM format and comprises $1353$ axial slices with $8.29 \rm GB$ size of memory in total. The range of grey values is from about $-4000$ to about $7500$. 

\begin{figure}
   \includegraphics[width = .7\textwidth]{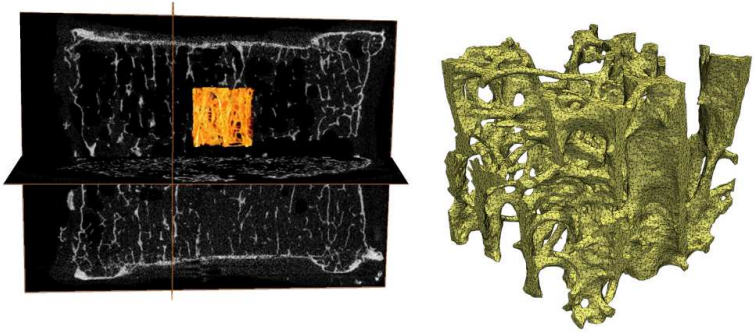}
 \caption{Multiplanar view of $\mu$CT data with the selected cubic sample in 3D-rendering (left) and STL-file of trabecular bone structure (right).}
 \label{fig4}
\end{figure}

Based on the DICOM data, the cubic sample of trabecular structure, depicted in Fig. \ref{fig4}, has been extracted using the strategy described in Section \ref{sec:ImageProc}. The cubic sample has about $5~{\rm mm}$ side length, comprises about $75,~000$ triangular faces and is intersection-free, fully manifold and with suitable element quality. Provided in STL-file format, the model is employed as input for the generation of the CFD-grid and for the underlying surface model for the volumetric FE-mesh, see Sections \ref{sec:CFD_Simu} and \ref{sec:FEM_Simu}, respectively.

\subsection{CFD-Simulation of the injection phase}
\label{sec:CFD_Simu}
	
\begin{figure}[ht]
  \centering \includegraphics[width=.35\textwidth]{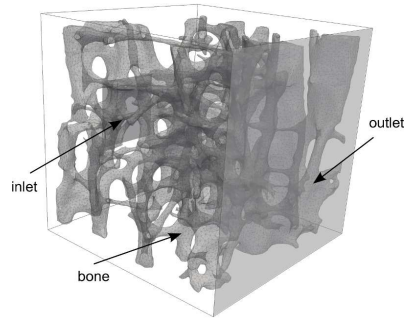}
 \caption{Boundary description of background mesh.}
 \label{fig5}
\end{figure}

The background mesh of the CFD-simulation has cubic dimensions with an edge length of $5.3~{\rm mm}$ and equally envelops the STL-geometry (see Fig. \ref{fig5}). The final mesh, where only space between the trabecular structure is meshed, has about $250,000$ cells. 

In order to define the boundary conditions of the CFD-model,  the notations given in Fig. \ref{fig5} are employed. A centred square on one of the vertical boundaries specifies the "inlet", where bone cement enters the computational domain. In our test case, we supposed that the bone cement injection starts three minutes after mixing the cement. The ambient temperature of bone cement is assumed to be $20~^\circ C$ prior to injection. Opposite to the inlet, bone marrow and bone cement can exit the computational domain through the "outlet". All other boundaries of the cube are denoted as "sides" and represent fixed walls where no fluid exchange is possible. The boundaries of the computational domain to the trabecular bone is referred to as "bone". Fluid flow through this boundary is also not permitted. Moreover, all boundaries except for the inlet remain adiabatic. The corresponding boundary conditions are summarised in Table \ref{tab2}.

Additionally to the boundary conditions, initial conditions have to be defined for the computational domain. To adopt the natural environment inside a human body, the initial temperature is set to $37^\circ{\rm C}$ and a pressure of $1~{\rm bar}$ is preset. Furthermore, it is assumed that the computational domain does not contain any bone cement in the initial state. Hence, the initial values for the phase fraction, the dissolution and the degree of cure are zero. Moreover, the initial velocity of bone marrow is set  $0~{\rm mm/s}$, thus no initial fluid flow is expected.

To complete the CFD model, transport properties have to be declared for both bone cement and bone marrow. To this end, the fluid model described in Section \ref{sec:fluid} is adopted for the bone cement phase. For simplicity and due to lack of data, bone marrow is modelled as a Newtonian fluid with constant viscosity. The corresponding material properties are listed in Table \ref{tab3} for both phases.

\begin{table}
\centering
  \caption{Boundary conditions of the CFD-model.}
  \label{tab2}
{\begin{tabular}{@{}p{30mm}@{\ \ \ \ \ }p{32mm}@{\ \ \ }p{28mm}@{\ \ \ }p{28mm}@{}}
 \hline 
                        & inlet                   & outlet             & sides/bone\\ 
 \hline \\[-3mm]
 pressure          & zero grad                     & fixed: $1\,\text{bar}$ & zero grad\\
 velocity            & fixed: $2.5\,\frac{\rm mm}{\rm s}$ & zero grad & fixed: $0\,\frac{\rm mm}{\rm s}$ \\ 
 phase fraction  & fixed: $1$                     & zero grad          & zero grad\\ 
 temperature     & fixed: $20\,^\circ\rm C$  & zero grad          & zero grad\\ 
 dissolution       & fixed: $156.7$               & zero grad          & zero grad\\ 
 degree of cure & fixed: $0.02$                & zero grad          & zero grad\\ 
 \hline \\
 \end{tabular}}
\end{table}

\begin{table}
\centering
  \caption{Material properties of the CFD-model.}
  \label{tab3}
{\begin{tabular}{@{}p{4cm}@{\ \ \ }p{3.1cm}@{\ \ \ }p{3.1cm}@{}}
\hline
 & bone cement & bone marrow\\ 
 \hline \\[-3mm]
kinematic viscosity 
   & see Section \ref{sec:fluid} & $4 \times 10^{-4}\,\text{m}^2/\text{s}$\\
density                  
   & $1.48\, \text{g}/\text{cm}^3$ &  $1.2\, \text{g}/\text{cm}^3$ \\
specific heat capacity
   & $1.2\, \text{kJ}/(\text{kg K})$ & $2\, \text{kJ}/(\text{kg K})$\\
thermal  conductivity  
   & $0.25\, \text{W}/(\text{m K})$ & $0.34\, \text{W}/(\text{m K})$\\
surface tension          
   & \multicolumn{2}{c}{$0.03\, \text{N}/\text{m}$}\\
\hline \\ \\
\end{tabular}}
\end{table}

\begin{figure}[ht]
  \centering
 \includegraphics[width=.99\textwidth]{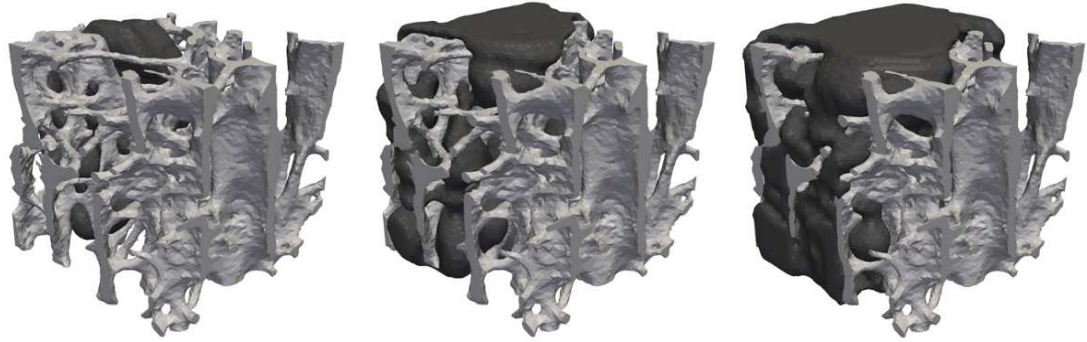}
 \caption{Bone cement distribution after 3s (left), 7s (middle) and 10s (right) of injection.}
 \label{fig6}
\end{figure}

The simulation of the case study was run on a workstation computer using parallel computing with a total number of four cores. Physical variables like pressure, velocity and temperature as well as the empirical variables phase fraction, degree of cure and progress of dissolution can be tracked throughout the whole simulation.  By way of example, the cement distributions after 3, 7 and 10 seconds of injection are shown in Fig. \ref{fig6}. 

Obviously, the cement penetrates greater cavities more easily than smaller ones. Regions, difficult to access, are hardly reached by bone cement. In our test case, the injection was supposed to stop at simulation time  $t_{CFD} = 10\,\text{s}$ and the results for temperature, cement distribution and degree of cure were passed to the FE-model.  

\subsection{Finite element simulation of bone cement curing}
\label{sec:FEM_Simu}

\begin{figure}[ht]
 \includegraphics[width=0.5\textwidth]{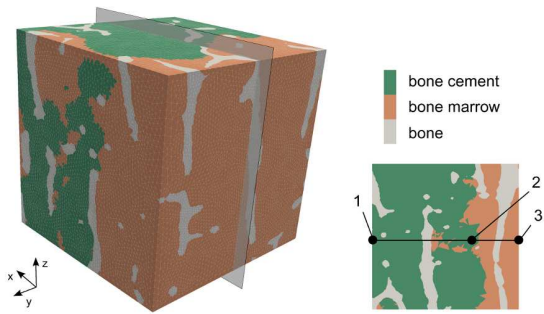}
 \caption{Finite element mesh (left) and intersecting plane (right): 1 - inlet, 2 - location near the bone cement interface, 3 - outlet.}
 \label{fig7}
\end{figure}

Subsequently to the injection simulation, the curing of bone cement inside the trabecular structure is investigated by finite element simulations. The FE-mesh of the case study again relies on the reconstructed computer model provided in STL-file format. An application of the meshing approach described in Section \ref{sec:FEM_modelling} results in a completely meshed cube, whereat the elements are divided into bone and pore elements. The resulting FE-mesh consists of approximately $770,000$ nodes and  $650,000$ tetrahedral elements. The bone cement distribution is transferred to the FE-mesh by application of the converter described in Section \ref{sec:FEM_modelling}. Thereby, the state provided by the CFD-simulation at simulation time of $10s$ is employed.  Furthermore, the distributions of temperature and degree of cure, which belong to the considered simulation time, are transferred to the model. By this, the initial conditions of bone cement and bone marrow are fully described. The resulting FE-mesh is divided into elements representing bone tissue ($\approx 210,000$ elements), bone marrow  ($\approx 160,000$ elements) and bone cement  ($\approx 280,000$ elements), respectively (see Fig. \ref{fig7}). The frame inserted in this figure indicates an intersecting plane, which is employed for subsequently visualisation of results. The highlighted positions represent the inlet, a location near the interface between bone cement and bone marrow as well as the outlet.

\begin{figure}[ht]
 \includegraphics[width=0.5\textwidth]{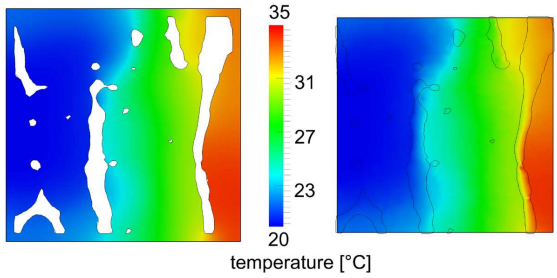}
 \caption{Comparison of temperature distribution between CFD- (left) and FE-model (right) at time of transition from CFD- to FE-simulation: contour plot of the temperature in the intersecting plane given in Fig. \ref{fig7} (right)).}
 \label{fig8}
\end{figure}

The first step of FE-simulation considers the temperature distribution of the whole mesh. This is due to the fact, that heat conduction in the bone tissue has not been computed within the CFD-simulation. Thus, firstly a steady-state heat transfer simulation is carried out, whereby the nodal temperatures of bone cement and bone marrow are fixed and only the steady-state heat conduction within the bone tissue is calculated. Fig. \ref{fig8} illustrates the resulting temperature distribution of this simulation step (right) as well as the corresponding temperature distribution provided by the CFD-simulation (left). A comparison of both temperature distributions indicates that the nodal temperatures of bone cement and bone marrow remain unchanged.

Having defined the initial conditions of the whole FE-model, the thermo-mechanically coupled process of bone cement curing is investigated. Towards this end, appropriate boundary conditions have to be defined. Here, two different sets of mechanical boundary conditions are investigated. In the first example, the displacements of nodes on the outer surfaces of the FE-mesh are assumed to be fixed perpendicular to their particular surfaces. However, the displacements in plane are not restricted such that nodes are allowed to move in plane of the outer surfaces. The second example equals the first one except that displacements on nodes of the outlet surface (cf. Fig. \ref{fig5}) are not restricted at all. Thus, this surface can deform without any restrictions. For both examples, an additional thermal boundary condition at the outlet (cf. Fig. \ref{fig5}) is responsible for the heat exchange to adjacencies. This boundary condition is adopted by heat convection with an ambient temperature $\theta_{amb} = 37 ^\circ C$ and a film coefficient $h_c = 1 {\rm mW/m^2\, K}$. All other surfaces of the cube remain adiabatic, thus, heat exchange is not allowed on those surfaces. 

Next, specific material behaviour is defined for the different materials of the FE-model. Firstly, the material model described in Section \ref{sec:solid} is assigned to the elements representing bone cement.  Moreover, simplified material behaviour is chosen for bone tissue and bone marrow. Both are modelled by an isotropic hyperelastic material using the Neo-Hookean ansatz
\begin{equation}
\label{eq:neo-hooke}
\rhotil \psi_{neo} = c_{10} \Big({\rm tr}(\C_G) - 3\Big) + \dfrac{9\,K}{2}\Big(\sqrt[3]{{\rm det}(\F)} - 1\Big)^2
\end{equation}
with $c_{10}$ and $K$ being a shear and bulk modulus, respectively. The operators "${\rm tr}$" and "${\rm det}$" denote the trace and the determinant of a second-rank tensor. The mechanical and thermal material parameters of bone tissue and bone marrow are summarised in Table \ref{tab4}. Specific material functions and corresponding parameters of the bone cement curing model are given in \cite{KolmederEtal2011}. 

\begin{table}[ht]
  \caption{Material properties for the FE-simulation.} 
  {\begin{tabular}{p{2.5cm} p{2.2cm} p{2.2cm} p{2.2cm} p{2.5cm} p{2.3cm}}
\hline
 & & & & & \\[-3mm]
     &  shear stiffness  
	 & bulk modulus 
	 & mass density 
	 & specific \ \ \ \ \ \ \ \ \ \ \ \ heat capacity 
	 & thermal \ \ \ \ \ \ \ \ \ \ \ \ conductivity\\[2mm]
\hline
& & & & & \\[-3mm]
bone
&  $1925 \ {\rm N/mm^2}$
&  $8333 \ {\rm N/mm^2}$
&  $2.0 \ \text{g}/\text{cm}^3$
&  $1300 \ \text{J}/(\text{kg K})$
&  $0.4 \ \text{W}/(\text{m K})$ \\[2mm]
bone marrow
&  $20 \ {\rm N/mm^2}$
&  $100 \ {\rm N/mm^2}$
&  $1.2\ \text{g}/\text{cm}^3$
&  $2000 \ \text{J}/(\text{kg K})$
&  $0.34\ \text{W}/(\text{m K})$ \\[2mm]
\hline
  \end{tabular}}
  \label{tab4}
\end{table}

The thermo-mechanically coupled FE-simulation of the curing process  was carried out on  a high performance computing cluster\footnote{Chemnitz High Performance Linux Cluster (CHiC), http://www.tu-chemnitz.de/chic/}, making use of the domain decomposition technology of MSC.MARC\textsuperscript{\textregistered}. In this case study, 8 compute nodes, each having 16 GB Ram and 2.64 GHz,were employed. The wall time for the curing simulation with a  total loadcase time of $\Delta t_{\rm FEM} = 625~{\rm s}$ was less than two hours. In the following, different exemplary results are provided. 

\begin{figure}[ht]
\begin{minipage}{0.48\textwidth}
\centering
\includegraphics[height=52mm]{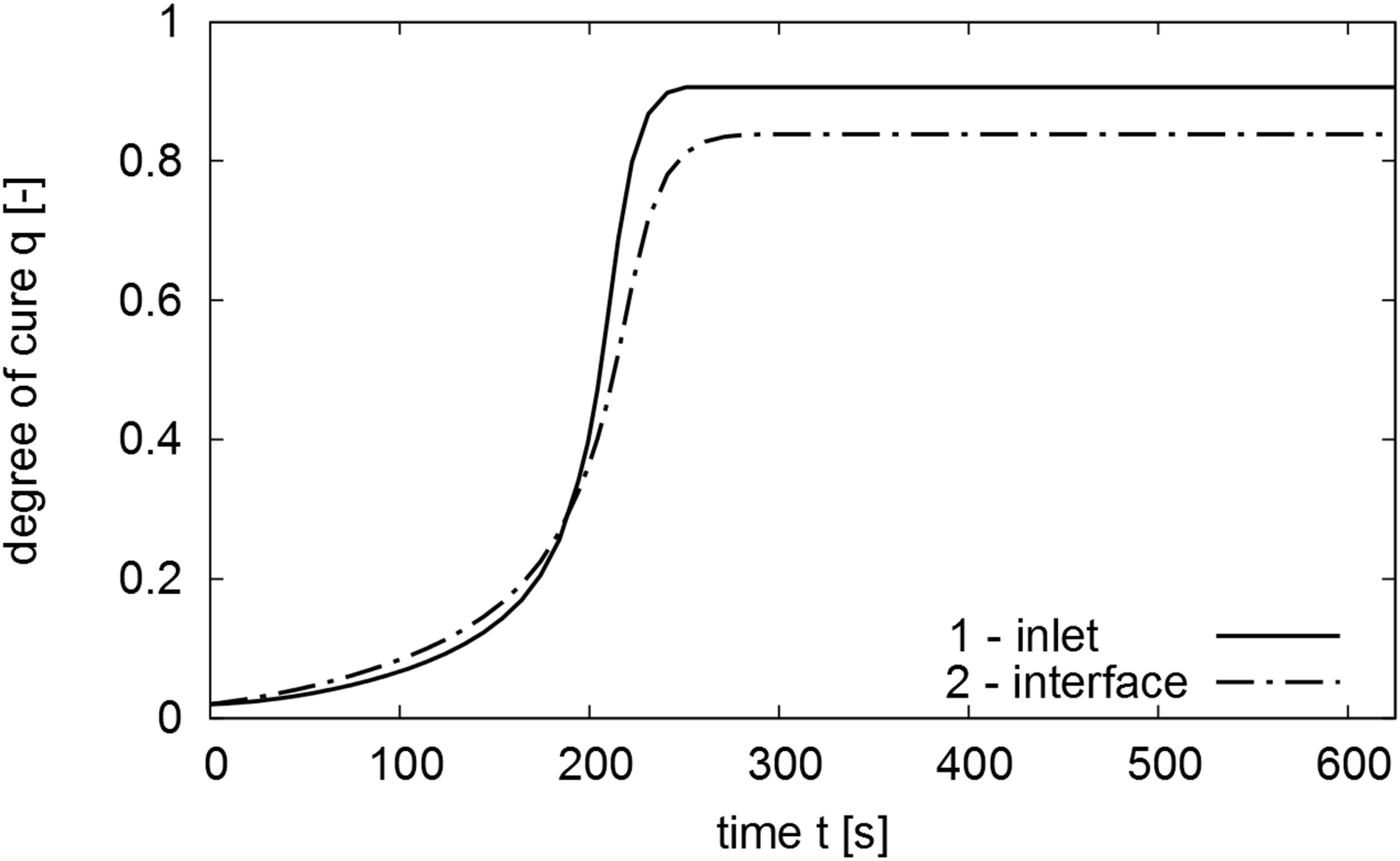}
\caption{Evolution of degree of cure at the inlet and bone cement/bone marrow interface  (cf. Fig. \ref{fig7}).}
\label{fig9}
\end{minipage}
\hfil
\begin{minipage}{0.45\textwidth}
\centering
\includegraphics[height=52mm]{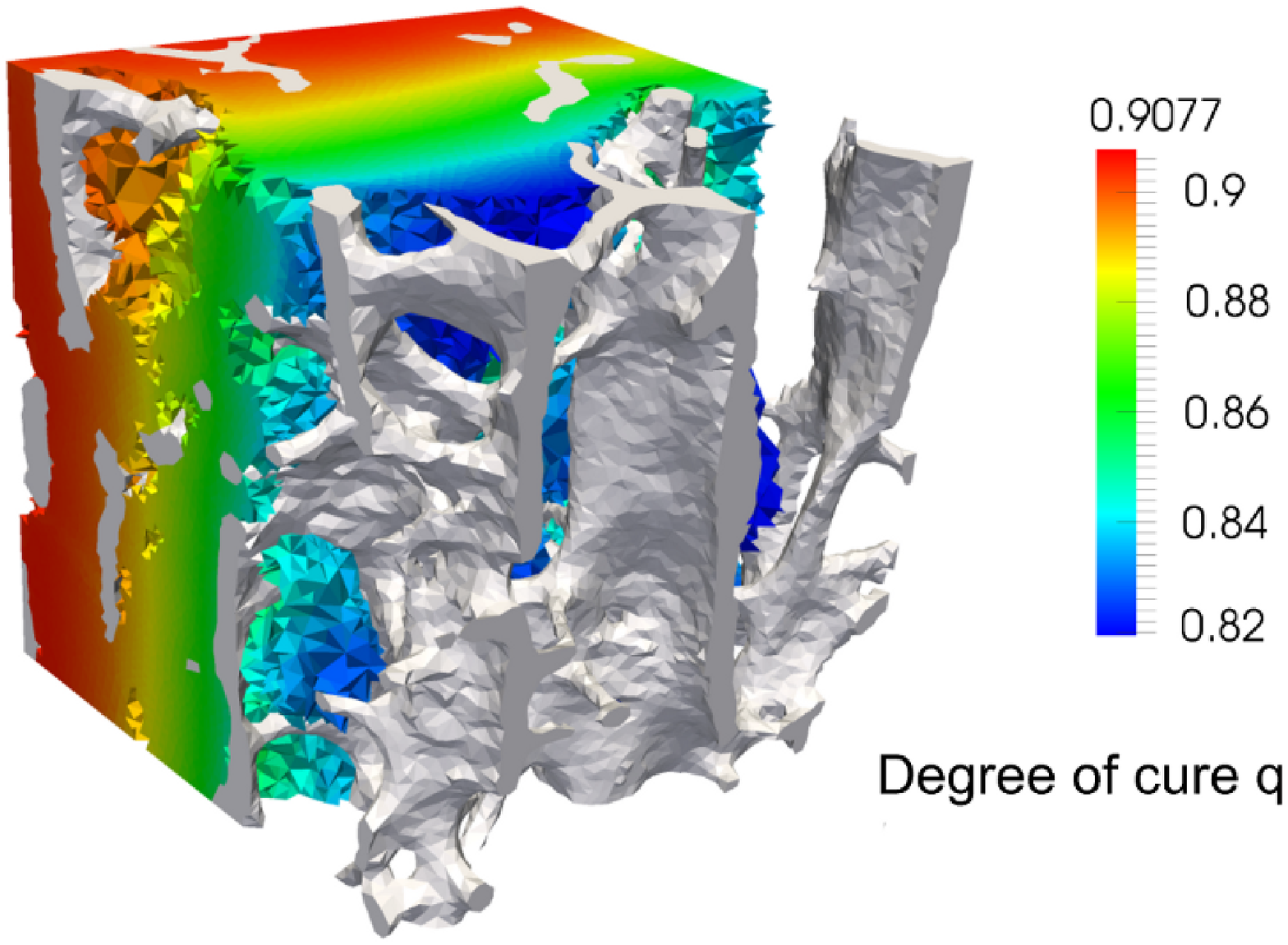}
\caption{Degree of cure after $625 s$ of finite element simulation.}
\label{fig10}
\end{minipage}
\end{figure}

First of all, the chemical process of polymerisation is evaluated. Polymerisation of acrylic bone cements, and thus curing of the material from a dough to a solid, is mainly governed by the degree of cure $q$ (cf. Eq. \eqref{eq:polymerisation}). The initial values for this variable were provided by the CFD-simulation described in Section \ref{sec:CFD_Simu}. Starting from this initial state, the degree of cure evolves until a certain saturated state is achieved. This behaviour is depicted in Fig. \ref{fig9}. The curves show the evolution of the degree of cure of bone cement at the inlet as well as at a location near the  bone cement/bone marrow interface  (cf. Fig. \ref{fig7}). Furthermore, the final distribution of the degree of cure at simulation time $t = 625\,s$ is depicted in Fig. \ref{fig10}.

\begin{figure}[ht]
 \includegraphics[width = 0.45\textwidth]{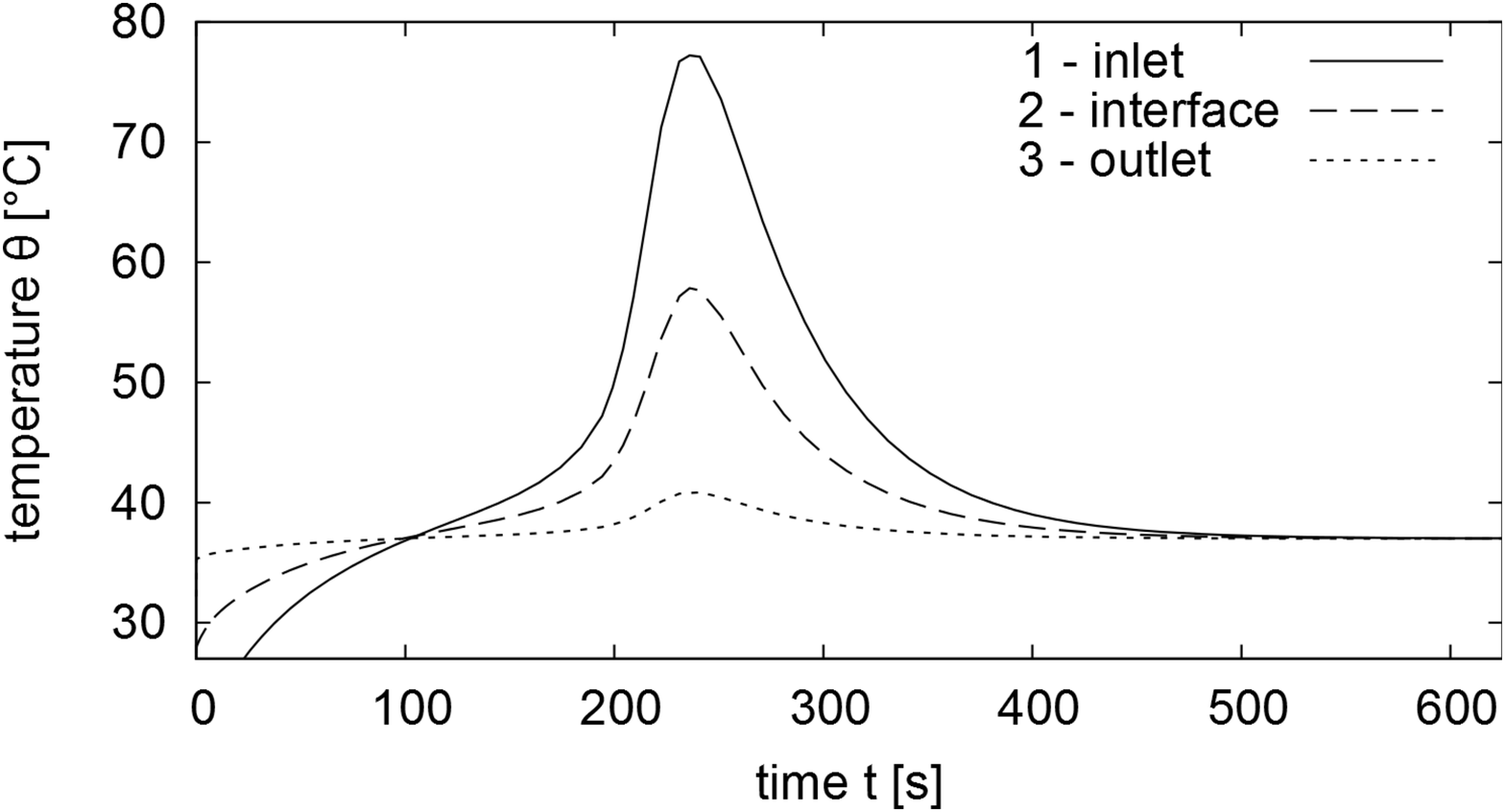}
 \caption{Evolution of temperature at different locations (cf. Fig. \ref{fig7}).}
 \label{fig11}
\end{figure}

Both Fig. \ref{fig9} and Fig. \ref{fig10} reveal the distinct temperature dependence of the degree of cure. On the one hand, a value of $q_{end} = 1$ is not achieved within this simulation. Consequently, the material has not been fully polymerised which results from the fact, that polymerisation takes place below the glass temperature ($T_g \approx 105 ^\circ C$). On the other hand, the degree of cure shows location dependent behaviour. This is due to different temperature evolutions at the specific locations. Local temperature evolutions in turn depend on the initial temperature distribution, the chosen thermal boundary conditions as well as the exothermic behaviour of the bone cement polymerisation. Those three main effects are covered by Fig. \ref{fig11}.

The initial temperature distribution of the FE-simulation is preset by the location dependent values obtained from the CFD-Simulation (cf. Fig. \ref{fig8}). Thereby, temperatures near the inlet have the lowest values. This is due to the definition that bone cement enters the computational domain with an ambient temperature of $\theta_{inlet} = 20 \, ^\circ C$. Within the first period of FE-simulation, heat convection leads the temperatures to rise up to body temperature. Subsequently, the impact of the exothermic chemical reaction can be observed. Thus, the temperature locally rises due to heat dissipation. The highest values ($\theta_{max} \approx 77 \, ^\circ C$) are obtained at the inlet, where no heat convection to the adjacencies is allowed. Due to the heat convection boundary condition at the outlet, temperature decreases until a homogeneous temperature field of $37 \, ^\circ C$ is achieved.  The general behaviour of the curing process and its impact on the thermal state correspond to the findings given in \cite{StanczykvanRietbergen2004}.

\begin{figure}[ht]
\begin{minipage}{0.47\textwidth}
\centering
  \includegraphics[width = \textwidth]{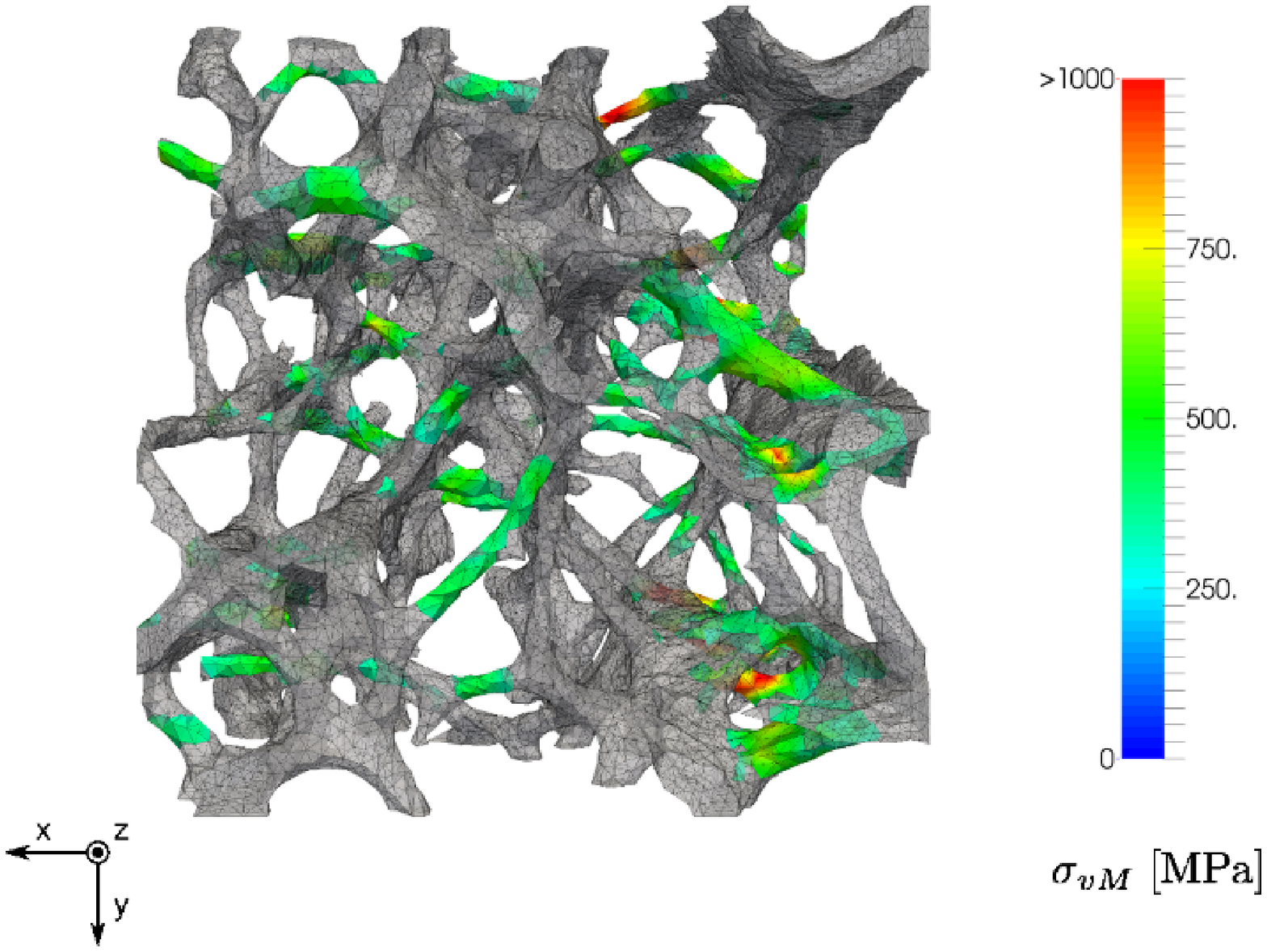}
 \caption{Equivalent von Mises stress of trabecular bone resulting from bone cement curing. All walls of the cube are fixed.}
 \label{fig12}
\end{minipage}
\hfil
\begin{minipage}{0.47\textwidth}
\centering
 \includegraphics[width = \textwidth]{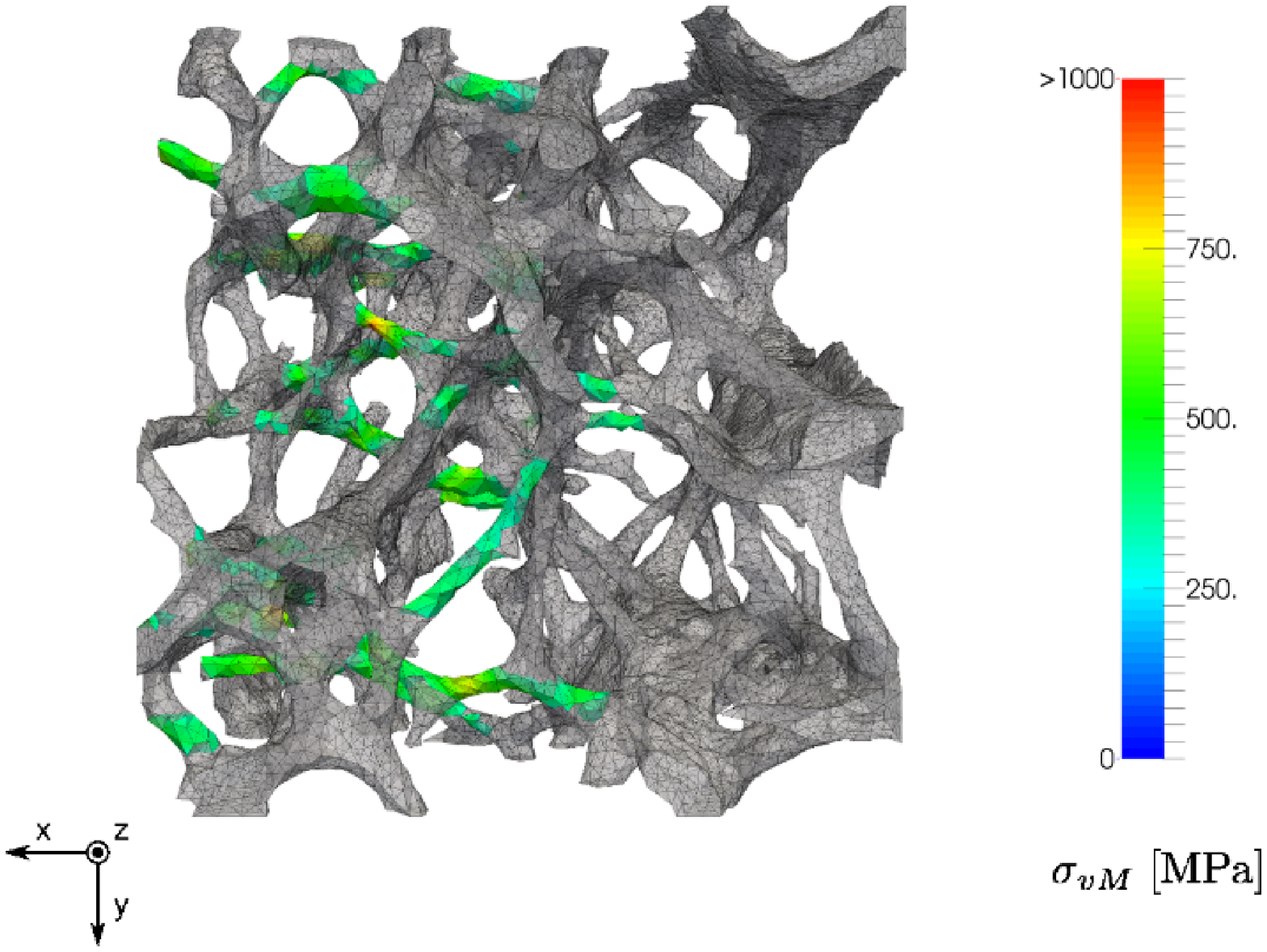}
 \caption{Equivalent von Mises stress of trabecular bone resulting from bone cement curing. Outlet wall (right surface) is not fixed.}
 \label{fig13}
\end{minipage}
\end{figure}

The curing of bone cement also has impact on the mechanics. Changes in temperature and the evolving degree of cure  lead to changes of the bone cement volume. Thereby, the chemical shrinkage due to polymerisation has a more dominant effect since it is known that acrylic bone cement can show a loss in volume of up $6-7 \%$ in stress-free state \cite{KuehnEtal2005}. Due to this local change of volume, stresses arise inside bone tissue. For the case of fully fixed surface walls, Fig. \ref{fig12} depicts the equivalent von Mises stresses that remain inside the bone tissue after $625~s$ of the curing simulation. Accordingly, Fig. \ref{fig13} illustrates the evolving stresses for the case of an unrestricted outlet surface. The viewing direction for both figures is in negative $z$-direction (cf. Fig. \ref{fig7}) and equivalent von Mises stresses with a magnitude of $\sigma_{vM}>300 \, \rm MPa$ are highlighted.

First of all it can be seen that stresses evolve in the bone tissue. This is only due to the bone cement's polymerisation process and the related chemical shrinkage of the material.  As can be seen from Fig. \ref{fig12}, which depicts the evolving stresses for the simulation with completely fixed outer surfaces, the highest values for the stresses occur in regions away from bone cement. This effect results from the chosen boundary conditions. Since nodes on the surfaces of the FE-model are fixed perpendicular to their surfaces, the total volume is forced to remain constant. Thus, the chemical shrinkage of bone cement results in expansion of the volume containing only bone tissue and bone marrow. Due to this expansion, bone tissue is stretched and stresses arise. However, if boundary conditions are chosen in a different way, the locations of trabecular failure vary. This behaviour is exemplified in Fig. \ref{fig13} for a simulation where the nodes on the outlet surface are not restricted. Hence, the total volume is allowed to change and high stress levels occur only in trabeculae completely surrounded by bone cement.

As a crucial result, it has to be noted that the stress levels in the bone tissue are non-physically high. Reported values for the yield stress of bone tissue are in the range of $100 - 200 \, \rm MPa$ (see, for example, \cite{BayraktarEtal2004,CezayirliogluaEtal1985,VerhulpEtal2008}). Due to the pure elastic material behaviour, which is assumed for bone tissue, these limits are locally exceeded. Hence, these locations point to failure of different trabeculae. More realistic stress predictions may be achieved by a refined material modelling, for example by including plasticity or damage to the bone material, or by the application of element birth and death technologies to account for fractures in the bone tissue. In this context, it was shown that even a simplified description of inelastic deformations in bone can lead to realistic apparent behaviour in microstructural finite element simulations \cite{vanRietbergenEtal1998,VerhulpEtal2008}. 

\begin{figure}[ht]
  \includegraphics[width=0.5\textwidth]{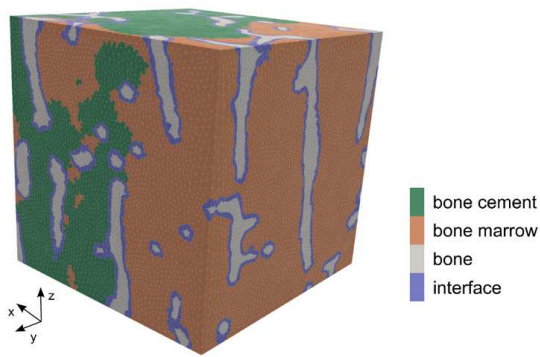}
 \caption{FE model incorporating an interface layer between bone and bone cement.}
 \label{fig14}
\end{figure}

\begin{figure}[ht]
\begin{minipage}{0.47\textwidth}
\centering
  \includegraphics[width = \textwidth]{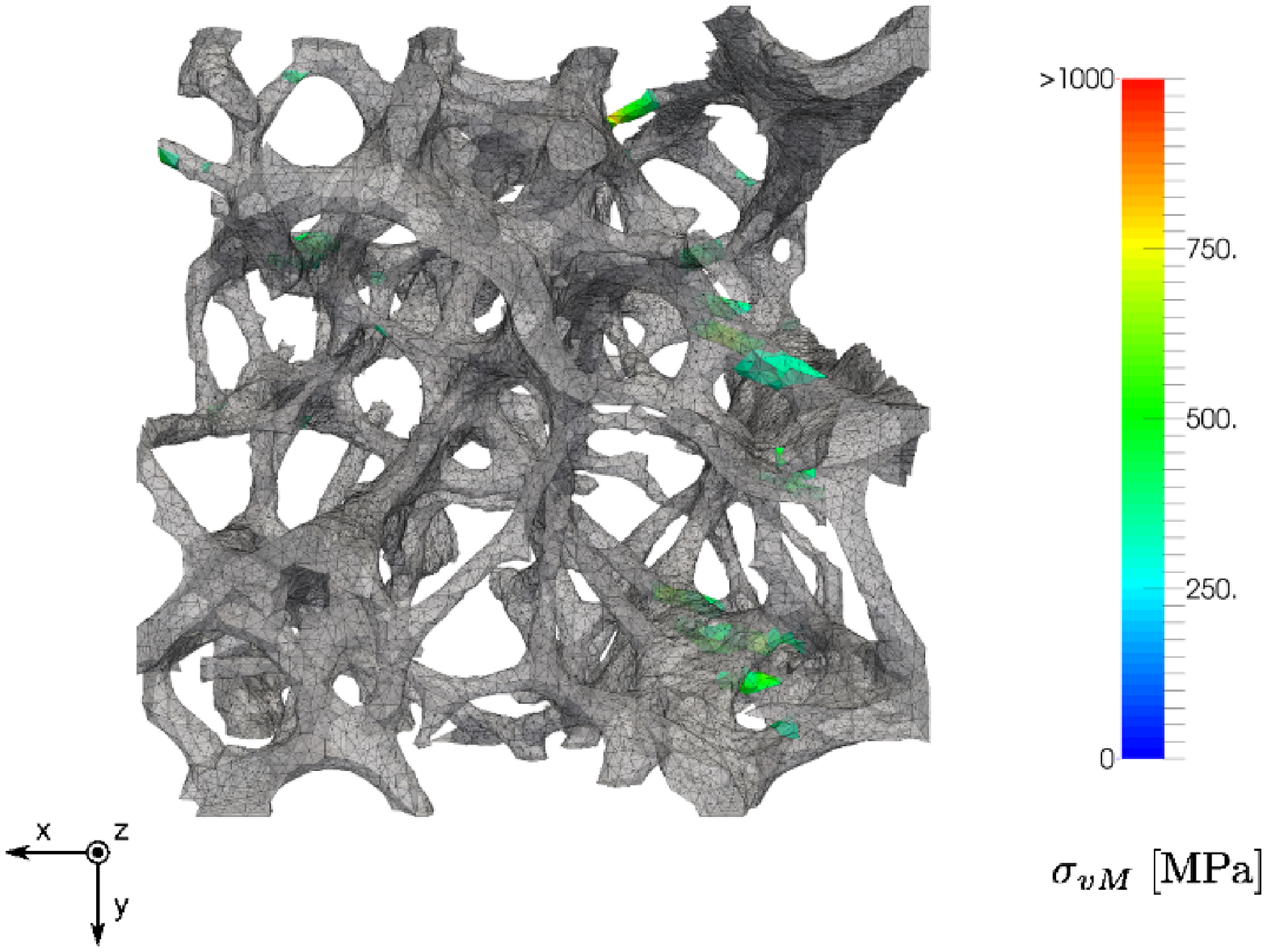}
 \caption{Equivalent von Mises stress of trabecular bone resulting from bone cement curing (soft interface layer between bone and bone cement, all walls of the cube are fixed).}
 \label{fig15}
\end{minipage}
\hfil
\begin{minipage}{0.47\textwidth}
\centering
 \includegraphics[width = \textwidth]{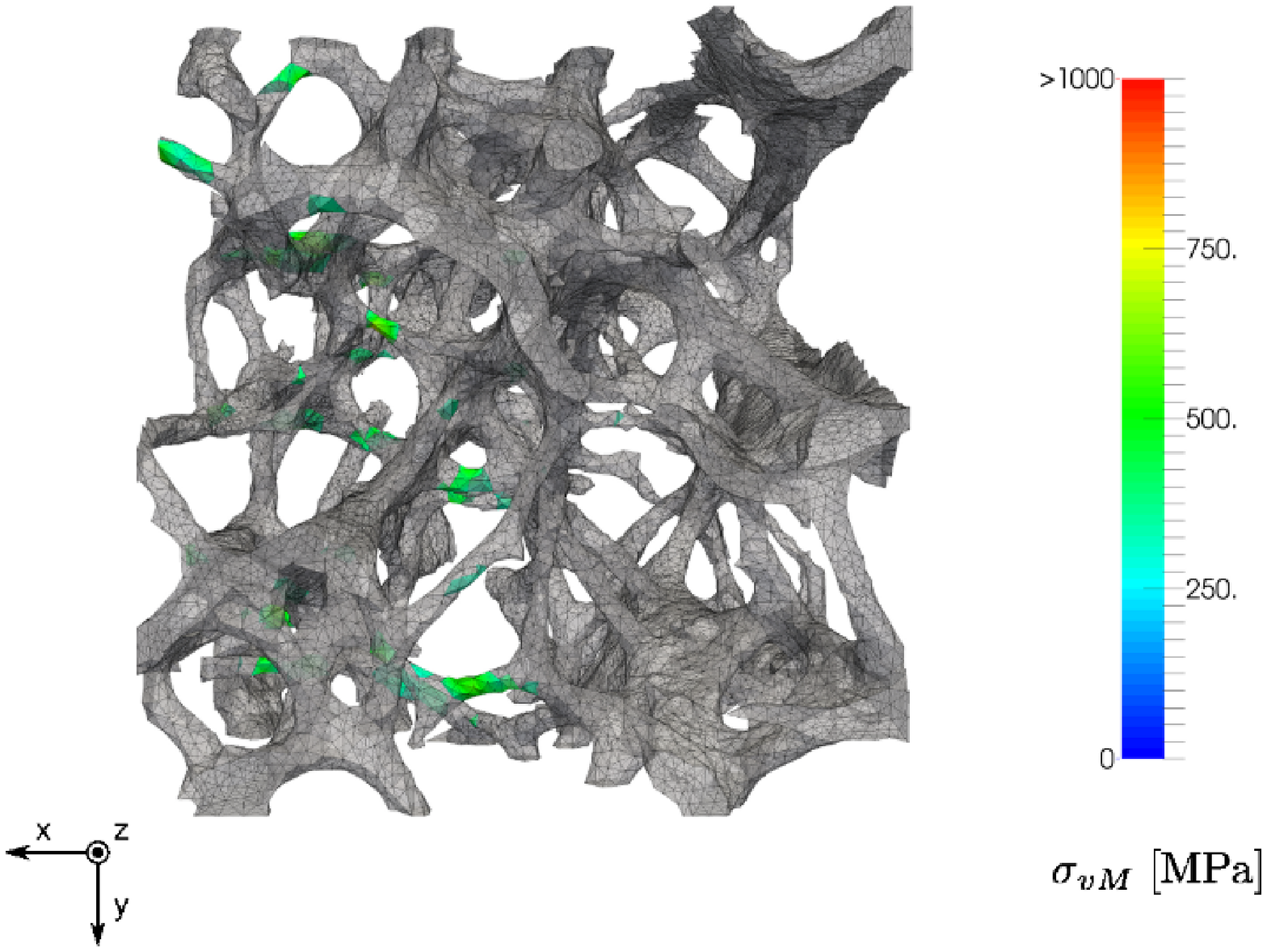}
 \caption{Equivalent von Mises stress of trabecular bone resulting from bone cement curing (soft interface layer between bone and bone cement, outlet wall is not fixed).}
 \label{fig16}
\end{minipage}
\end{figure}

Another potential cause of these high stress levels shall be investigated in more detail. As reported in the literature, bone cement and bone may not form a perfect bonding due to the hydrophobic behaviour of bone cement and the hydrophilic behaviour of bone and bone marrow \cite{KinzlEtal2012,MummeEtal2007}. Instead of a chemical bonding, the interface may only be formed by mechanical contact and even have pores. This clearly reduces the load transmission between both materials and the bone/bone cement compound is only connected by mechanical interlocking \cite{KinzlEtal2012}. In order to investigate the impact of a non-perfect interface bonding, a modified model has been generated where an interface layer between bone cement and bone is included (see Fig. \ref{fig14}). It is assumed that the interface layer behaves like bone marrow. Thus, the interface layer has a comparably low stiffness. With this modified model, the two simulations with completely fixed outer surfaces and one free surface, respectively, were repeated, and the evolving stresses in the bone tissue are depicted in Figs. \ref{fig15} and \ref{fig16}.

Compared to the results of the original model with perfect bonding behaviour between bone and bone cement (cf. Figs. \ref{fig12} and \ref{fig13}), both figures reveal that lower stress values and fewer affected trabeculae occur in the modified model. However, some regions are still exposed to stress levels exceeding the critical limit of $100 -200 \rm \, MPa$ and thus risk of fractures inside the cancellous bone structure is still predicted. Even if the impact of the interface layer could be shown in these investigations, it has to be noted that the in vivo behaviour is affected by the bone cement composition and the resulting ability of bone cement to attach to the bone structure. Further numerical studies on this issue can be found in \cite{KinzlEtal2012}.

\section{Concluding remarks}
  \label{sec:closing}

In this paper, a computational framework for the simulation of bone cement injection and curing processes related to the surgical treatment of vertebroplasty was presented. The framework includes the generation of 3D-computer models  based on patient specific data obtained from $\mu$CT imaging, thermo-mechanically coupled CFD-simulations of bone cement injection into trabecular structures and thermo-mechanically coupled FE-simulations of bone cement curing. Through all simulations stages, a detailed material description of the acrylic bone cement's material behaviour was used. A non-linear and process-depending fluid flow model was employed to represent the initially liquid state of the biomaterial. The curing process was described by a non-linear viscoelastic model including cure kinetics, temperature dependencies and evolving mechanical parameters, among others. 

The capabilities of the presented approach were demonstrated with the help of a test case made of an osteoporotic cancellous bone sample. Here, different results from injection and curing simulations were presented. The CFD-simulation of the injection process revealed reasonable filling patterns and the results referring to the thermo-chemical curing process correspond to findings reported in the literature. Moreover, it was shown that chemical shrinkage of the bone cement leads to residual stresses in the bone tissue and, furthermore, can induce damage inside the cancellous bone structure. The distinct impact of different boundary conditions and an assumed bone/bone cement interface layer on these results was demonstrated as well.

In our future work it is intended to refine the modelling approach and to perform further studies regarding the injection, curing and long term behaviour of acrylic bone cements. Aspects of model refinement are, for example, the investigation of different sets of boundary conditions with an even more practical relevance and an enhanced material description for the bone tissue, which includes inelastic or damage effects. Besides different studies concerning the impact of varied process parameters, the modelling and simulation strategy will also be employed to estimate effective anisotropic physical properties of the cancellous bone structure and the bone/bone cement compound. This includes the identification of effective permeability, heat conduction and stiffness properties.

\section*{Acknowledgements}

The authors kindly acknowledge the financial support by the German Research Foundation (DFG) within the project PAK 273. Furthermore, we thank Werner Schm\"olz for providing the $\mu$CT data and Thomas R. Blattert for fruitful discussions.


\end{document}